\documentclass[letterpaper, 11pt]{article}

\usepackage{amsmath,amssymb,amsthm}
\usepackage{color}
\usepackage{enumitem}
\usepackage[left=2.5cm, right=2.5cm, top=2.5cm, bottom = 2.5cm]{geometry}
\usepackage{graphicx}
\usepackage[numbers,sort]{natbib}
\usepackage[hidelinks]{hyperref}
\usepackage[center]{subfigure}

\definecolor{lightGreen}{rgb}{0.6, 1, 0.6}
\definecolor{darkGreen}{rgb}{0.0, 0.5, 0.0}
\definecolor{lightBlue}{rgb}{0.8, 0.8, 1}
\definecolor{lightRed}{rgb}{1, 0.8, 0.8}

\newcommand{\red}[1]{\textcolor{red}{#1}}
\newcommand{\blue}[1]{\textcolor{blue}{#1}}
\newcommand{\green}[1]{\textcolor{darkGreen}{#1}}

\newcommand{\bolda}{\boldsymbol{a}}
\newcommand{\boldf}{\boldsymbol{f}}
\newcommand{\boldg}{\boldsymbol{g}}
\newcommand{\boldx}{\boldsymbol{x}}
\newcommand{\boldy}{\boldsymbol{y}}

\newcommand{\boldv}{\boldsymbol{v}}
\newcommand{\boldw}{\boldsymbol{w}}

\newcommand{\boldnu}{\boldsymbol{\nu}}

\newcommand{\boldxi}{\boldsymbol{\xi}}
\newcommand{\boldeta}{\boldsymbol{\eta}}

\newcommand{\intomega}{\int_{\Omega}}
\newcommand{\intomegah}{\int_{\Omega_h}}
\newcommand{\intgamma}{\int_{\Gamma}}
\newcommand{\intgammah}{\int_{\Gamma_h}}

\newcommand{\nablagamma}{\nabla_{\Gamma}}
\newcommand{\nablagammah}{\nabla_{\Gamma_h}}

\DeclareMathOperator{\volume}{volume}
\DeclareMathOperator{\area}{area}
\DeclareMathOperator{\diam}{diam}

\DeclareMathOperator{\Tr}{Tr}

\newtheorem{definition}{Definition}
\newtheorem{lemma}{Lemma}
\newtheorem{proposition}{Proposition}
\newtheorem{remark}{Remark}
\newtheorem{theorem}{Theorem}

\begin{document}

\title{Virtual element method for elliptic bulk-surface PDEs in three space dimensions}

\author{%
{\sc
Massimo Frittelli\thanks{Corresponding author. Email: massimo.frittelli@unisalento.it}} \\[2pt]
Department of Mathematics and Physics ``E. De Giorgi'', University of Salento\\
Via per Arnesano, 73100 Lecce, Italy \\[6pt]
{\sc Anotida Madzvamuse}\thanks{Email: A.Madzvamuse@sussex.ac.uk}\\[2pt]
Department of Mathematics, School of Mathematical and Physical Sciences,\\
University of Sussex, Brighton, BN1 9QH, UK
{\sc and}\\[6pt]
{\sc
Ivonne Sgura\thanks{Corresponding author. Email: ivonne.sgura@unisalento.it}} \\[2pt]
Department of Mathematics and Physics ``E. De Giorgi'', University of Salento\\
Via per Arnesano, 73100 Lecce, Italy 
}

\maketitle

\begin{abstract}
In this work we present a novel bulk-surface virtual element method (BSVEM) for the numerical approximation of elliptic bulk-surface partial differential equations (BSPDEs) in three space dimensions.  The BSVEM is based on the discretisation of the bulk domain into polyhedral elements with arbitrarily many faces. The polyhedral approximation of the bulk induces a polygonal approximation of the surface.  Firstly, we present a geometric error analysis of bulk-surface polyhedral meshes independent of the numerical method. Then, we show that BSVEM has optimal second-order convergence in space, provided the exact solution is $H^{2+3/4}$ in the bulk and $H^2$ on the surface, where the additional $\frac{3}{4}$ is due to the combined effect of surface curvature and polyhedral elements close to the boundary.  We show that general polyhedra can be exploited to reduce the computational time of the matrix assembly.  To demonstrate optimal convergence results, a numerical example is presented on the unit sphere.
\end{abstract}

\section*{Keywords}
Bulk-surface PDEs;  Polyhedral meshes;  Bulk-surface virtual
element method; Convergence.

\section*{Mathematics Subject Classification}
 65N12, 65N15, 65N30, 65N50

\section{Introduction}
In this work we introduce the bulk-surface virtual element method (BSVEM) for the numerical approximation of elliptic bulk-suface partial differential equations (BSPDEs) in three space dimensions of the following form:
\begin{equation}
\label{elliptic_problem}
\begin{cases}
- \Delta u (\boldx) + u(\boldx) = f(\boldx), \qquad \boldx \in \Omega;\\
- \Delta_\Gamma v (\boldx) + v(\boldx) + \nabla u(\boldx)\cdot \boldnu(\boldx)  = g(\boldx), \qquad \boldx \in \Gamma;\\
\hspace*{2.5mm} \nabla u(\boldx)\cdot \boldnu(\boldx) = -\alpha u(\boldx) + \beta v(\boldx), \qquad \boldx \in \Gamma,\\
\end{cases}
\end{equation}
where $\Omega \subset\mathbb{R}^3$ is an open set such that $\Gamma = \partial \Omega$ is a smooth surface,  $\Delta$ is the Laplace operator in $\Omega$, $\Delta_{\Gamma}$ is the Laplace-Beltrami operator on $\Gamma$,  $\boldnu$ is the outward unit normal vector field on $\Gamma$,  $\alpha,  \beta > 0$ and $f:\Omega\rightarrow\mathbb{R}$ and $g: \Gamma\rightarrow\mathbb{R}$ are data.  The problem \eqref{elliptic_problem} is taken from \cite{elliottranner2013finite} and is the prototype of \emph{coupled bulk-surface partial differential equations} (BSPDEs),  a class of problems that is recently drawing attention in the literature.  More generally, given a number $d\in\mathbb{N}$ of space dimensions, a system of BSPDEs comprises of $m\in\mathbb{N}$ PDEs posed in the \emph{bulk} $\Omega \subset\mathbb{R}^d$, coupled with $n\in\mathbb{N}$ PDEs posed on the surface $\Gamma := \partial \Omega$ through either linear or non-linear coupling, see for instance \cite{Madzvamuse_2016}. The quickly growing interest toward BSPDEs arises from the numerous applications of such PDE problems in different areas, such as cellular biological systems (\cite{Mackenzie_2016, Elliott_2017, cusseddu2019, Paquin_Lefebvre_2019}), fluid dynamics (\cite{Bianco_2013, Lee_2015, burman2016cut}),  and plant biology (\cite{ryder2019sensor}) among many other applications.

Among the various state-of-the art numerical methods for the spatial discretisation of BSPDEs existing in the literature we mention bulk-surface finite elements (BSFEM) (\cite{elliottranner2013finite, Madzvamuse_2015, Madzvamuse_2016, Kovacs_2016}), trace finite elements (\cite{gross2015trace}), cut finite elements (\cite{burman2016cut}), discontinuous Galerkin methods (\cite{chernyshenko2018hybrid}),  kernel collocation method (\cite{chen2019kernel}),  and closest point method (\cite{macdonald2013simple}). 

The purpose of the present paper is to introduce a novel \emph{bulk-surface virtual element method} (BSVEM) for the spatial discretisation of elliptic BSPDEs in $d=3$ space dimensions.  The BSVEM is a substantial extension of the recently introduced virtual element method (VEM) for the numerical approximation of several classes of partial differential equations on flat domains (\cite{beirao2013basic}) or surfaces (\cite{frittelli2018virtual}). The key feature of VEM is that of being a \emph{polyhedral} method, i.e.  it handles elements of a quite general polyhedral shape, rather than just of tetrahedral shape (\cite{beirao2013basic}). The success of virtual elements is due to several advantages arising from polyhedral mesh generality, such as: (i) computationally cheap mesh pasting (\cite{benkemoun2012anisotropic, chen2014memory, frittelli2018virtual}), (ii) efficient adaptive algorithms (\cite{cangiani2014adaptive}),  (iii) flexible approximation of the domain and its boundary (\cite{dai2007n}),  (iv) nonconforming elements (\cite{gardini2019nonconforming}), and (iv) the possibility of enforcing higher regularity to the numerical solution (\cite{antonietti2016c, da2013arbitrary, brezzi2013virtual}). Thanks to these advantages, several extensions of the original VEM for the Poisson equation (\cite{beirao2013basic}) were developed for numerous PDE problems, such as heat (\cite{vacca2015virtual}) and wave equations (\cite{vacca2016virtual}), reaction-diffusion systems (\cite{adak2019convergence}), Cahn-Hilliard equation (\cite{antonietti2016c}), Stokes equation (\cite{da_veiga_2017_stokes}), Helmholtz equation (\cite{MASCOTTO2019445}), linear elasticity (\cite{da2013virtual}), plate bending (\cite{brezzi2013virtual}), fracture problems with geophysical applications (\cite{benedetto2016globally,  fumagalli2021performances}), eigenvalue problems (\cite{mora2015virtual}) and many more.

On one hand,  our proposed numerical methodology combines the VEM for the bulk equations (\cite{da2017high}) with the surface virtual element method (SVEM) (\cite{frittelli2018virtual}) for the surface equations.  On the other hand, the numerical method extends the two-dimensional BSVEM introduced in \cite{frittelli2021bulk}.  A marked difference with the work presented in \cite{frittelli2021bulk} is that the surface PDEs were solved using the (one-dimensional) surface finite elements, in this work, we employ virtual elements for both bulk and surface PDEs.  Here, the method relies on an arbitrary polyhedral discretisation of the bulk and its corresponding induced polygonal approximation of the surface.  To the best of our knowledge, this kind of geometrical approximation is novel in the literature.  In the special case of tetrahedral meshes, the method boils down to the BSFEM (\cite{elliottranner2013finite, Madzvamuse_2015}). 

The theoretical novelty of the present study is threefold. 
Firstly,  we provide a geometric error analysis of polyhedral bulk-surface meshes that is independent of the numerical method and applies, in principle, to any method based on polyhedral bulk-surface meshes.
Secondly, we carry out a full error analysis of the BSVEM.  The proposed method possesses optimal second-order convergence provided the numerical solution is $H^{2+3/4}(\Omega)$ in the bulk instead of the usual requirement of $H^2(\Omega)$, see \cite{elliottranner2013finite}.  However,  our analysis requires such extra regularity only in the simultaneous presence of a curved boundary $\Gamma$ and non-tetrahedral elements close to the boundary, a novel case.  We point out that such extra regularity comes for free in most models and applications, where the domains are smooth and the solutions are infinitely differentiable.  Whether such higher regularity is also necessary, it remains an open problem. 
Thirdly,  in the case $\beta = 0$ in \eqref{elliptic_problem}, the first equation in \eqref{elliptic_problem} becomes a bulk-only PDE with non-zero Neumann conditions.  Hence, a by-product of the proposed analysis is that the lowest-order VEM for bulk-only elliptic PDEs in 3D retains optimal convergence in the simultaneous presence of a curved boundary $\Gamma$, non-tetrahedral elements close to $\Gamma$, and non-zero Neumann data.  Interestingly enough, this problem was fully addressed only in specific cases.  For example, in the simplest case of tetrahedral meshes (FEM),  the result was proven in \cite{Bartels_2004}.  In the case of general polyhedral meshes (VEM), the seminal work (\cite{da2017high}) is confined to polyhedral domains.  Then, in \cite{beirao2019curved} and \cite{dassi2021bend} the authors consider a VEM in 2D with curved edges and a VEM in 3D with curved faces, respectively,  to take out the geometric error.  In \cite{bertoluzza2019} the authors introduce a 2D VEM with suitable algebraic corrections that account for curved boundaries.  The present work finally addresses the 3D case and does not require any geometric or algebraic correction of the VEM. 

In addition, we show that the usage of suitable polyhedra drastically reduces the computational time of matrix assembly on equal meshsize in comparison to the tetrahedral BSFEM.  This property, which already holds true in the 2D case (\cite{frittelli2021bulk}), is even more accentuated in 3D.  Similar results are obtained in the literature through other methods, such as trace FEM (\cite{gross2015trace}) or cut FEM (\cite{burman2016cut}). 

The structure of our paper is as follows.
In Section \ref{sec:bspdes} we derive the weak formulation of problem \eqref{elliptic_problem} and we state existence, uniqueness and regularity results.
In Section \ref{sec:geometric_analysis} we introduce polyhedral bulk-surface meshes and analyse the geometric error.  
In Section \ref{sec:bsvem_method} we introduce the BSVEM for problem \eqref{elliptic_problem}. In Section \ref{sec:convergence_analysis} we carry out the convergence analysis.
In Section \ref{sec:mesh_advantage} we show that polyhedral meshes can  significantly reduce the computational time of the matrix assembly.
In Section \ref{sec:example_elliptic_bs_sphere} we provide a numerical example on the sphere to demonstrate the optimal convergence.
In Appendix A we provide basic definitions and results required for the analysis.

\section{Weak formulation, existence and regularity}
\label{sec:bspdes}

To obtain the weak formulation of  \eqref{elliptic_problem}, we multiply the first two equations of \eqref{elliptic_problem} by two test functions $\alpha\varphi \in H^1(\Omega)$ and $\beta\psi \in H^1(\Gamma)$, respectively,  then we apply Green's formula in the bulk $\Omega$ and Green's formula on the curved manifold $\Gamma$ (\cite{dziuk2013finite}). We obtain the following formulation: find $u \in H^1(\Omega)$ and $v \in H^1(\Gamma)$ such that
\begin{equation}
\label{elliptic_problem_weak_form_intermediate}
\begin{cases}
\vspace{2mm}
\displaystyle \alpha\intomega \Big(\nabla u \cdot \nabla \varphi + u\varphi\Big) = \alpha\intomega f\varphi + \alpha\intgamma \dfrac{\partial u}{\partial\boldnu}\varphi;\\
\displaystyle \beta\intgamma \Big(\nablagamma v\cdot\nablagamma \psi + v\psi\Big) + \beta\intgamma \dfrac{\partial u}{\partial \boldnu} \psi = \beta\intgamma g\psi,
\end{cases}
\end{equation}
for all $\varphi \in H^1(\Omega)$ and $\psi \in H^1(\Gamma)$. By using the third equation of \eqref{elliptic_problem} in \eqref{elliptic_problem_weak_form_intermediate} and summing over the equations, we obtain the following weak formulation: find $(u,v)\in H^1(\Omega) \times H^1(\Gamma)$ such that 
\begin{equation}
\label{elliptic_problem_weak_form}
b((u,v); (\varphi,\psi)) = \intomega f\varphi + \intgamma g\psi,
\end{equation}
for all $(\varphi,\psi) \in H^1(\Omega) \times H^1(\Gamma)$, where $b((u,v); (\varphi,\psi)): (H^1(\Omega) \times H^1(\Gamma))^2 \rightarrow\mathbb{R}$ is the bilinear form defined by
\begin{equation*}
b((u,v); (\varphi,\psi)) = \displaystyle \alpha\intomega \Big(\nabla u \cdot \nabla \varphi + u\varphi\Big) + \beta\intgamma \Big(\nablagamma v\cdot\nablagamma \psi + v\psi\Big) + \intgamma (\alpha u - \beta v)(\alpha\varphi-\beta\psi).
\end{equation*}
The variational formulation \eqref{elliptic_problem_weak_form} fulfils the following result on existence, uniqueness and regularity found in  \cite{elliottranner2013finite}.
\begin{theorem}[Existence, uniqueness and regularity (\cite{elliottranner2013finite})]
If $\Gamma$ is a $\mathcal{C}^3$ surface, $f\in L^2(\Omega)$ e $g\in L^2(\Gamma)$,  the variational problem \eqref{elliptic_problem_weak_form} has a unique solution $(u,v) \in H^2(\Omega) \times H^2(\Gamma)$ that fulfils the following bound
\begin{align}
\label{exact_solution_bound}
\|(u,v)\|_{H^2(\Omega)\times H^2(\Gamma)} \leq C \|(f,g)\|_{L^2(\Omega) \times L^2(\Gamma)} .
\end{align}
\end{theorem}
Thanks to elliptic regularity, it is also possible to show that if $f\in H^1(\Omega)$ and $g\in H^1(\Omega)$, the regularity improves to
\begin{align}
\label{exact_solution_bound_H^3}
\|(u,v)\|_{H^3(\Omega)\times H^2(\Gamma)} \leq C \|(f,g)\|_{L^2(\Omega) \times L^2(\Gamma)} .
\end{align}

\section{Geometric analysis}
\label{sec:geometric_analysis}
In this section we introduce bulk-surface polyhedral meshes and we analyse the geometric approximation error.  The present analysis is independent of the numerical method and applies, in principle, to any polyhedral method for BSPDEs.

\subsection{Polyhedral bulk-surface meshes}

\begin{figure}
\begin{center}
\subfigure[Illustration of the \emph{bulk} $\Omega$ enclosed by the \emph{surface} $\red{\Gamma}$, the \emph{narrow band} $\green{U_\delta}$ and the Fermi stripe $\blue{U}$.]{\includegraphics[scale=.6]{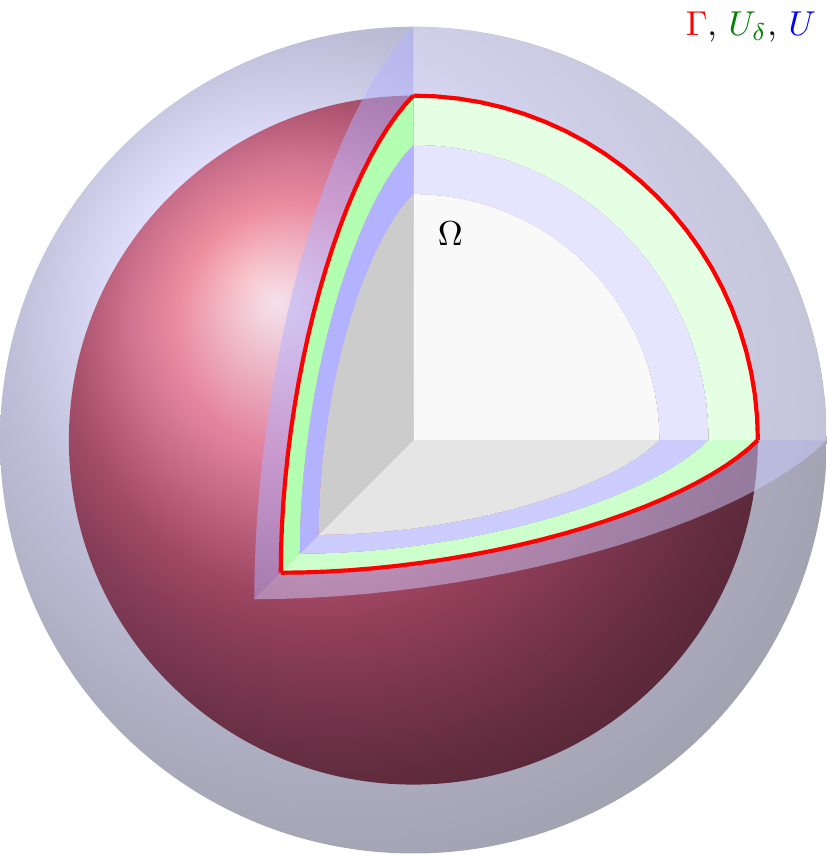}\label{fig:exact_domain}}
\hspace*{1cm}
\subfigure[Illustration of the \emph{discrete bulk} $\Omega_h$ enclosed by the \emph{discrete surface} $\red{\Gamma_h}$, the \emph{discrete narrow band} $\green{\Omega_B}$ and the Fermi stripe $\blue{U}$.]{\includegraphics[scale=.6]{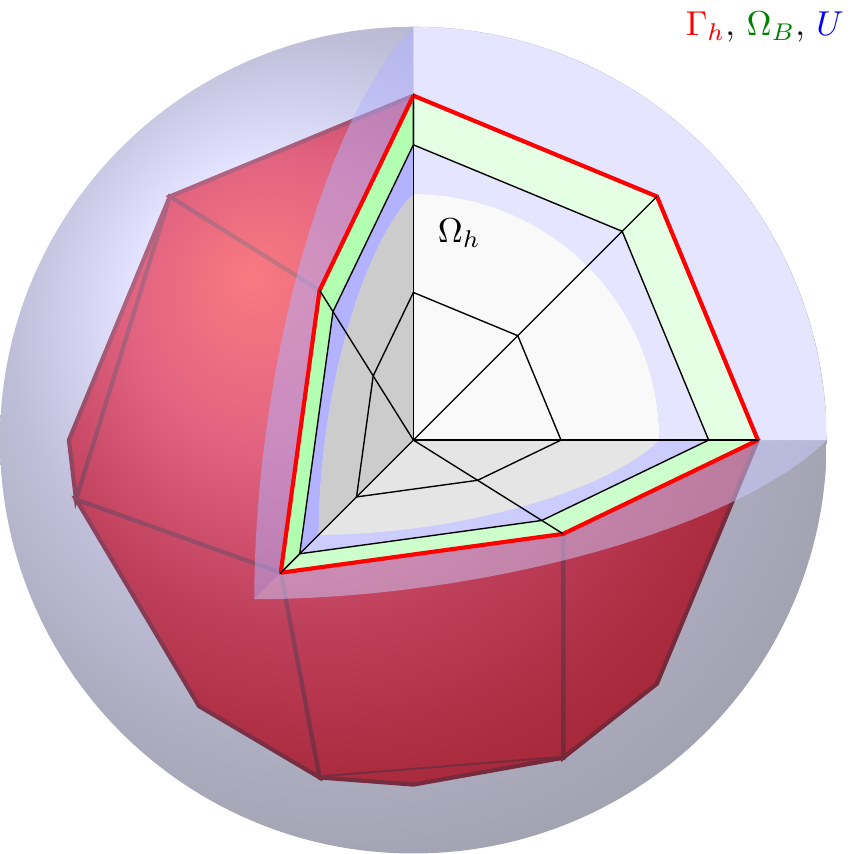}\label{fig:discrete_domain}}
\end{center}
\caption{Illustration of the continuous domain, the discrete domain and related notations.}
\label{fig:mesh_illustration}
\end{figure}

Let $h>0$ be a positive number called \emph{meshsize} and let $\Omega_h = \cup_{E\in\mathcal{E}_h} E$ be a polyhedral approximation of the bulk $\Omega$, where $\mathcal{E}_h$ is a set of non-degenerate compact polyhedra. The polyhedral bulk $\Omega_h$ automatically induces a polygonal approximation $\Gamma_h$ of $\Gamma$, defined by $\Gamma_h = \partial \Omega_h$, exactly as in the case of tetrahedral meshes,  see \cite{elliottranner2013finite}. Notice that we can write $\Gamma_h = \cup_{F \in \mathcal{F}_h} F$, where $\mathcal{F}_h$ is the set of the faces of $\Omega_h$ that constitute $\Gamma_h$. We assume that:
\begin{enumerate}[label=\textnormal{(F\arabic*)}]
\item the diameter of each element $E\in\mathcal{E}_h$ does not exceed $h$;\label{B1}
\item for any two distinct elements or faces, their intersection is either empty, or a common vertex, or a common edge, or a common face.
\item all nodes of $\Gamma_h$ lie on $\Gamma$;
\item every face $F\in\mathcal{F}_h$ is contained in the Fermi stripe $U$ of $\Gamma$ (see Fig.  \ref{fig:mesh_illustration}).\label{B4}
\end{enumerate}
\begin{enumerate}[label=\textnormal{(V\arabic*)}]
\item there exists $\gamma_1 > 0$ such that every $E \in \mathcal{E}_h$ and every face $F$ of $E$ is star-shaped with respect to a ball (with center $x_E$ and $x_F$ respectively) of radius $\gamma_1 h_E$ and $\gamma_1 h_F$ respectively, where $h_E$ and $h_F$ are the diameters of $E$ and $F$, respectively;\label{A1} 
\item there exists $\gamma_2 > 0$ such that for all $E \in \mathcal{E}_h$ and for and every face $F$ of $E$, the distance between any two nodes of $E$ or $F$ is at least $\gamma_2 h_E$ or $\gamma_2 h_F$, respectively.\label{A2}
\end{enumerate}
Assumptions (F1)-(F4) are standard in the SFEM literature, see for instance \cite{dziuk2013finite}, while assumptions (V1)-(V2) are standard in the VEM literature, see for instance \cite{beirao2013basic}. The combined assumptions (F1)-(V2) will prove sufficient in our bulk-surface setting. In the following definitions and results we provide the necessary theory for estimating the geometric error arising from the boundary approximation.
\begin{definition}[Essentials of polyhedral bulk-surface meshes]
\label{def:boundary_elements}
An element $E\in\mathcal{E}_h$ is called an \emph{exterior element} if it has at least a face or an edge contained in $\Gamma_h$, otherwise $E$ is called an \emph{interior element}.
Let $\Omega_B$ be the \emph{discrete narrow band} defined as the union of the exterior elements of $\Omega_h$ as illustrated in Fig. \ref{fig:discrete_domain}. From Assumption (F4), for any face $F$ contained in $\Gamma_h$ we have that $\bolda(F) \subset \Gamma$, where $\bolda$ is the normal projection defined in \eqref{lmm:fermi}. 
\end{definition}

Observe that, for sufficiently small $h>0$, the discrete narrow band $\Omega_B$ is contained in the Fermi stripe $U$ as shown in Fig.  \ref{fig:discrete_domain}. Let $N \in\mathbb{N}$ and let $\boldx_i$, $i=1,\dots,N$, be the nodes of $\Omega_h$. Let $M \in\mathbb{N}$, $M < N$ and assume that the nodes of $\Gamma_h$ are $\boldx_k$, $k=1,\dots,M$, i.e. the first $M$ nodes of $\Omega_h$. Throughout the paper we need the following \emph{reduction matrix} $R \in\mathbb{R}^{N\times M}$ defined as $R := [I_M; 0]$, where $I_M$ is the $M\times M$ identity matrix. The reduction matrix $R$ fulfils the following two properties:
\begin{itemize}
\item For $\boldv \in \mathbb{R}^{N}$, $R^T\boldv \in \mathbb{R}^M$ is the vector with the first $M$ entries of $\boldv$;
\item For $\boldw \in \mathbb{R}^M$, $R\boldw \in \mathbb{R}^N$ is the vector whose first $M$ entries are those of $\boldw$ and the other $N-M$ entries are $0$.
\end{itemize}
In what follows, we will use the matrix $R$ for an optimised implementation of the BSVEM.

\subsection{Variational crime}
\begin{figure}
\begin{center}
\subfigure[Every boundary face of the element $\green{E}$ is contained in the Fermi stripe $\blue{U}$ of $\red{\Gamma}$.]{\label{fig:proof_step_1}\includegraphics[scale=1.2]{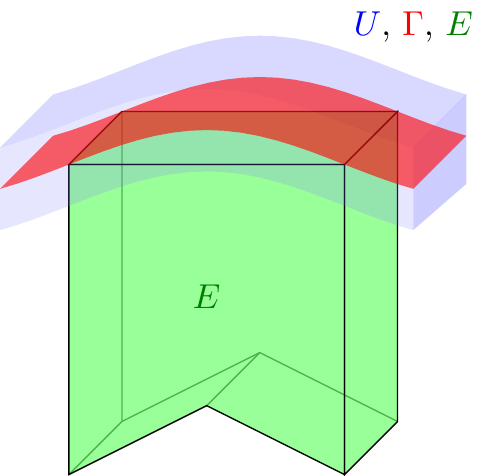}}
\subfigure[The element $\green{E}$ is split into tetrahedra, like $\red{T}$ shown in the picture, having $\boldx_E$ as a vertex. ]{\label{fig:proof_step_2}\includegraphics[scale=1.2]{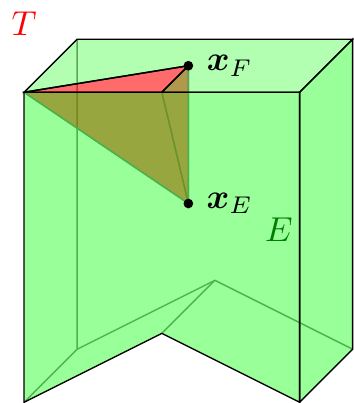}}
\end{center}
\begin{center}
\subfigure[The tetrahedron $\red{T}$ is replaced by its exact counterpart $\red{\breve{T}}$.]{\label{fig:proof_step_3}\includegraphics[scale=1.2]{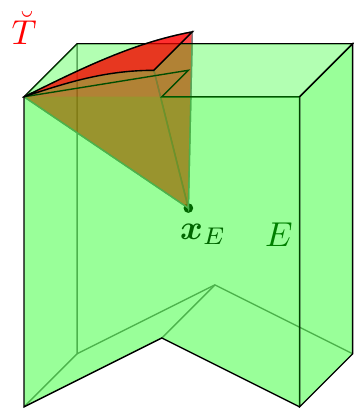}}
\hspace*{10mm}
\subfigure[By repeating for each $T$, the exact element $\breve{E}$ is obtained.]{\label{fig:proof_step_4}\includegraphics[scale=1.2]{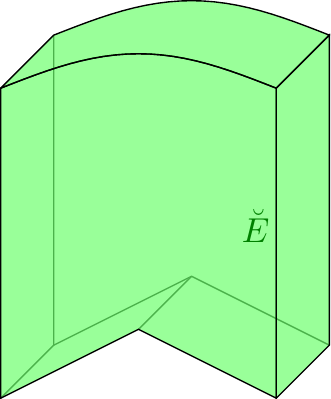}}
\end{center}
\caption{Steps of the construction of the exact element $E$ corresponding to a given element $E$, following Lemma  \ref{lmm:homeomorphism}.  The symbols $\blue{U}, \red{\Gamma}, \green{E}$ are colour-matched with the figure.}
\label{fig:proof}
\end{figure}

We now consider the geometric error due to the boundary approximation. Since the surface variational crime in surface virtual elements is well-understood (\cite{frittelli2018virtual}), we will mainly focus on the variational crime in the bulk. To this end, it is useful to analyse the relation between any element $E\in\mathcal{E}_h$ and a suitably defined \emph{exact element} $\breve{E}$ (a curved version of $E$),  see Fig. \ref{fig:proof} for an illustration. For the special case of tetrahedral meshes with at most one boundary face per element,  $\breve{E}$ is rigorously defined, there exists a diffeomorphism between $E$ and $\breve{E}$ and this diffeomorphism is linearly close to the identity with respect to the meshsize, see \cite{elliottranner2013finite}.  In the more general case when $E$ has more than four faces and/or multiple boundary faces,  we will show the existence of a mapping between $E$ and a suitably defined $\breve{E}$ with slightly weaker regularity, which is sufficient for our purposes.

\begin{lemma}[Domain parametrisation]
\label{lmm:homeomorphism}
Let $\mathcal{E}_h$ fulfil assumptions (F1)-(V2). There exists a homeomorphism $G:\Omega_h \rightarrow \Omega$ such that $G \in W^{1,\infty}(\Omega_h)$ and
\begin{align}
\label{estimate_1}
&G|_{\Gamma_h} = \bolda|_{\Gamma_h};\\
\label{estimate_2}
&G|_{\Omega_h \setminus \Omega_B} = Id;\\
\label{estimate_3}
&\|JG - Id\|_{L^\infty(\Omega_B)} \leq Ch;\\
\label{estimate_4}
&\|\det(JG) - 1\|_{L^\infty(\Omega_B)} \leq Ch;\\
\label{estimate_5}
&\|G - Id\|_{L^\infty(\Omega_B)} \leq Ch^2,
\end{align}
where $\bolda$ is the normal projection defined in Lemma \ref{lmm:fermi}, $JG$ is the Jacobian of $G$ and $C$ is a constant that depends on $\Gamma$ and the constants $\gamma_1$, $\gamma_2$ are those considered in Assumptions (V1)-(V2).  Even if restricted to a single element $E\in\mathcal{E}_h$, $G$ might not be a diffeomorphism unless $E$ is a tetrahedron.

\begin{proof}
Consider a bulk element $E\in\mathcal{E}_h$ and assume that all of the faces of $E$ that are contained in $\Gamma_h$ are also in the Fermi stripe $U$, see Fig.  \ref{fig:proof_step_1}. 
Pick a face $F$ of $E$ and let $x_E$ and $x_F$ be as in Assumption (V1). By joining $x_E$ and $x_F$ with the midpoints of two consecutive edges of $F$, a tetrahedron $T$ is obtained,  see Fig.  \ref{fig:proof_step_2}.
By proceeding in this fashion,  $E$ can be subdivided into a finite amount $N_E$ of tetrahedra $T_1,\dots,T_{N_E}$ that are quasi-uniform thanks to the geometric assumptions (V1)-(V2). 
Then, replace each $T_i$ by its exact (curved) counterpart $\breve{T}_i$ as defined in \cite{elliottranner2013finite}, see Fig.  \ref{fig:proof_step_3}. 
The exact element $\breve{E}$ is then defined by replacing each $T_i$ by its curved counterpart $\breve{T}_i$, see Fig. \ref{fig:proof_step_4}. 
The claimed map $G$ is constructed piecewise by applying \cite[Proposition 4.7]{elliottranner2013finite} for all the $T_i$'s of each $E\in\mathcal{E}_h$.
If restricted to a single $T_i$, the map $G$ is a diffeomorphism \cite[Proposition 4.7]{elliottranner2013finite}.
\end{proof}
\end{lemma}

Thanks to Lemma \ref{lmm:homeomorphism} it is possible to define bulk- and surface-lifting operators.
\begin{definition}[Bulk- and surface-lifting operators]
Given $V:\Omega_h\rightarrow \mathbb{R}$ and $W:\Gamma_h \rightarrow\mathbb{R}$, their \emph{lifts} are defined by $V^\ell := V \circ G^{-1}$ and $W^\ell := W \circ G^{-1}$, respectively. Conversely, given $v:\Omega \rightarrow\mathbb{R}$ and $w:\Gamma \rightarrow\mathbb{R}$, their \emph{inverse lifts} are defined by $v^{-\ell} := v \circ G$ and $w^{-\ell} := w \circ G$,
respectively, with $G:\Omega_h \rightarrow\Omega$ being the mapping defined in Lemma \ref{lmm:homeomorphism}.
\end{definition}
Lemma \ref{lmm:homeomorphism} also enables us to show the equivalence of Sobolev norms under lifting as illustrated next.
\begin{lemma}[Equivalence of norms under lifting]
There exists two constants $c_2 > c_1 > 0$ depending on $\Gamma$ and $\gamma_2$ such that, for all $V: \Omega_h \rightarrow \mathbb{R}$ and for all $W: \Gamma_h \rightarrow \mathbb{R}$,
\begin{align}
\label{equivalence_omega_l2}
&c_1 \|V^\ell\|_{L^2(\Omega)} \hspace*{3mm} \leq \|V\|_{L^2(\Omega_h)}  \hspace*{3mm} \leq c_2\|V^\ell\|_{L^2(\Omega)};\\
\label{equivalence_omega_h1}
&c_1 |V^\ell|_{H^1(\Omega)} \leq |V|_{H^1(\Omega_h)} \leq c_2|V^\ell|_{H^1(\Omega)};\\
\label{equivalence_gamma_l2}
&c_1 \|W^\ell\|_{L^2(\Gamma)}  \hspace*{2mm} \leq \|W\|_{L^2(\Gamma_h)}  \hspace*{2mm} \leq c_2\|W^\ell\|_{L^2(\Gamma)};\\
\label{equivalence_gamma_h1}
&c_1 |W^\ell|_{H^1(\Gamma)} \leq |W|_{H^1(\Gamma_h)} \leq c_2|W^\ell|_{H^1(\Gamma)};\\
\label{equivalence_gamma_h2}
& \hspace{25mm} |W|_{H^2(\Gamma_h)} \leq c_2|W^\ell|_{H^2(\Gamma)} + c_2 h |W^\ell|_{H^1(\Gamma)}.
\end{align}
\begin{proof}
Estimates \eqref{equivalence_omega_l2}-\eqref{equivalence_omega_h1} follow by using the map $G$ introduced in Lemma \ref{lmm:homeomorphism} in the proof of \cite[Proposition 4.9]{elliottranner2013finite}. A proof of \eqref{equivalence_gamma_l2}-\eqref{equivalence_gamma_h2} is in \cite[Lemma 4.2]{dziuk2013finite}.
\end{proof}
\end{lemma}

We are ready to estimate the effect of lifting on bulk- and surface integrals.
\begin{lemma}[Geometric error of lifting]
\label{lmm:geometric_errors_bilinear forms}
If $u, \varphi\in H^{1}(\Omega)$, then
\begin{align}
\label{lifting_error_bulk_nabla}
&\left |\intomega \nabla u \cdot \nabla \varphi - \intomegah \nabla u^{-\ell} \cdot \nabla \varphi^{-\ell} \right| \leq Ch|u|_{H^1(\Omega_B^\ell)}|\varphi|_{H^1(\Omega_B^\ell)},\\
\label{lifting_error_bulk_mass}
&\left |\intomega u \varphi - \intomegah u^{-\ell} \varphi^{-\ell} \right| \leq Ch \|u\|_{L^{2}(\Omega_B^\ell)}\|\varphi\|_{L^{2}(\Omega_B^\ell)},
\end{align} 
where $C$ depends on $\Gamma$, $\gamma_1$ and $\gamma_2$. If $v,\psi\in H^1(\Gamma)$, then
\begin{align}
\label{lifting_error_surf_nabla}
&\left|\intgamma \nablagamma v \cdot \nablagamma \psi - \intgammah \nablagammah v^{-\ell} \cdot \nablagammah \psi^{-\ell} \right| \leq Ch^2|v|_{H^1(\Gamma)} |\psi|_{H^1(\Gamma)};\\
\label{lifting_error_surf_mass}
&\left|\intgamma v \psi - \intgammah v^{-\ell} \psi^{-\ell} \right| \leq Ch^2\|v\|_{L^2(\Gamma)} \|\psi\|_{L^2(\Gamma)},
\end{align}
where $C$ depends on $\Gamma$, $\gamma_1$ and $\gamma_2$.

\begin{proof}
To prove \eqref{lifting_error_bulk_nabla}-\eqref{lifting_error_bulk_mass} it suffices to use the bulk geometric estimates \eqref{estimate_2}-\eqref{estimate_4} in the proof of \cite[Lemma 6.2]{elliottranner2013finite}. A proof of \eqref{lifting_error_surf_nabla}-\eqref{lifting_error_surf_mass} can be found in \cite{dziuk2013finite}.
\end{proof}
\end{lemma}

\begin{remark}[Polyhedral meshes and curved boundaries]
\label{rmk:lifting_regularity}
From Lemma \ref{lmm:homeomorphism} we know that the mapping $G$ might not be a diffeomorphism in the simultaneous presence of general polyhedral elements and curved boundaries.  This issue does not arise in the absence of curved boundaries (\cite{da2017high}), when $G$ is the identity by construction,  or in the absence of non-tetrahedral elements (\cite{Kovacs_2016}). This implies that, in the simultaneous presence of curved boundaries and non-tetrahedral elements, the lifting operator does not preserve the Sobolev regularity of functions. That is to say,  for $E\in\mathcal{E}_h$ the inverse lift of an $H^2(\breve{E})$ function is not, in general, $H^2(E)$. Now, since our analysis requires \emph{full regularity of the exact solution mapped on the polyhedral domain}, we need an alternative mapping instead of the lifting.  Hence, we consider the \emph{Sobolev estension}.
\end{remark}

\begin{lemma}[Geometric error of Sobolev extension]
\label{lmm:geometric_error_of_sobolev_extension}
There exist $C>0$ such that
\begin{align}
\label{error_between_lift_and_extension_l2}
&\|\tilde{u} - u^{-\ell}\|_{L^2(\Omega_h)} \leq Ch^2\|u\|_{H^{2+1/4}(\Omega)}, \quad \forall\ u\in H^{2+1/4}(\Omega);\\
\label{error_between_lift_and_extension}
&|\tilde{u} - u^{-\ell}|_{H^1(\Omega_h)} \leq Ch^{\frac{3}{2}}\|u\|_{H^2(\Omega)} + Ch\|u\|_{H^{2+3/4}(\Omega)}, \quad \forall\ u\in H^{2 +1/4}(\Omega).
\end{align}
\begin{proof}
By using \eqref{sobolev_extension}, \eqref{sobolev_inequality_holder} with $\gamma = \frac{3}{4}$, \eqref{estimate_2} and \eqref{estimate_5} we have that
\begin{align}
\label{my_lemma_proof_0}
\|\tilde{u}-u^{-\ell}\|_{L^2(\Omega_h)} = \|\tilde{u} - \tilde{u}\circ G\|_{L^2(\Omega_h)} \leq C\|\tilde{u}\|_{H^{3/2+3/4}(\Omega_h)}\|(Id-G)^{3/4}\|_{L^2(\Omega_h)} \notag\\
= C\|u\|_{H^{3/2+3/4}(\Omega)}\|(Id-G)^{3/4}\|_{L^2(\Omega_B)} \leq C\|u\|_{H^{3/2+3/4}(\Omega)} |\Omega_B|^{1/2}\|Id-G\|_{L^\infty(\Omega_B)}^{3/4} \notag\\
\leq Ch^{\frac{1}{2}} h^{\frac{3}{2}}\|u\|_{H^{3/2+3/4}(\Omega)} = Ch^2\|u\|_{H^{2+1/4}(\Omega)},
\end{align}
which proves \eqref{error_between_lift_and_extension_l2}. Notice that, in the last line of \eqref{my_lemma_proof_0}, the $h^\frac{1}{2}$ term is the effect of the Sobolev extension being exact except on the discrete narrow band $\Omega_B$. Using \eqref{narrow_band_inequality}, \eqref{sobolev_extension}, \eqref{estimate_3} and \eqref{estimate_5} we have that
\begin{align}
\label{my_lemma_proof_1}
|\tilde{u} - u^{-\ell}|_{H^1(\Omega_h)} = \|\nabla \tilde{u} - (JG^T\ \nabla\tilde{u})\circ G\|_{L^2(\Omega_h)} \leq\\
\notag
\|(Id-JG^T\circ G)\|_{L^\infty(\Omega_h)}\|\nabla \tilde{u}\|_{L^2(\Omega_B)} + \|JG^T\circ G\|_{L^\infty(\Omega_h)}\|\nabla\tilde{u} - \nabla\tilde{u}\circ G\|_{L^2(\Omega_h)} \leq\\
\notag
Ch\|\nabla\tilde{u}\|_{L^2(\Omega_B)} + C\|\nabla\tilde{u} - \nabla\tilde{u}\circ G\|_{L^2(\Omega_h)} \leq Ch^{\frac{3}{2}}\|u\|_{H^2(\Omega)} + C\|\nabla\tilde{u} - \nabla\tilde{u}\circ G\|_{L^2(\Omega_h)}.
\end{align}
Since $\tilde{u} \in H^{2+1/2+\gamma}(\Omega_h)$, then $\nabla\tilde{u} \in H^{1+1/2+\gamma}(\Omega_h)$. Hence, by reasoning as in \eqref{my_lemma_proof_0} we have that
\begin{equation}
\label{my_lemma_proof_3}
\|\nabla\tilde{u} - \nabla\tilde{u}\circ G\|_{L^2(\Omega_h)} \leq C_\gamma h^{\frac{1}{2} + 2\gamma}\|u\|_{H^{2+1/2+\gamma}(\Omega)}.
\end{equation}
By substituting \eqref{my_lemma_proof_3} into \eqref{my_lemma_proof_1} we get the desired estimate.
\end{proof}
\end{lemma}

\section{The Bulk-Surface Virtual Element Method (BSVEM)}
\label{sec:bsvem_method}
In this section we introduce the Bulk-Surface Virtual Element Method (BSVEM) for problem \eqref{elliptic_problem}. 

\subsection{Virtual element space on polygons and polyhedra}
We start by defining virtual element spaces on polygons and polyhedra by following \cite{da2017high}, but we simplify the presentation,  as the present work is confined to first-degree elements.  We start from the two dimensional spaces. Let $F$ be a polygon that, without loss of generality,  lies in $\mathbb{R}^2$. A preliminary virtual element space on $F$ is given by
\begin{equation}
\tilde{\mathbb{V}}(F) := \Big\{v \in H^1(F) \cap \mathcal{C}^0(F) \Big|\ v_{|e} \in\mathbb{P}_1(e),\ \forall\ e\in\text{edges}(F) \land \Delta v \in \mathbb{P}_1(F) \Big\},
\end{equation}
where $\mathbb{P}_1(F)$ is the space of linear polynomials on the polygon $F$.  Let us consider the elliptic projection $\Pi^\nabla_F : \tilde{\mathbb{V}}(F) \rightarrow \mathbb{P}_1(F)$ defined by
\begin{equation}
\vspace{2mm}
\displaystyle\int_F \nabla (v - \Pi^\nabla_F v) \cdot \nabla p_1 = 0\qquad \forall\ p_1 \in \mathbb{P}_1(F) \quad  \land \quad \displaystyle\int_{\partial F} (v - \Pi^\nabla_F v) = 0.
\end{equation}
Thanks to Green's formula, the operator $\Pi_F^\nabla$ is computable, see \cite{ahmad2013equivalent} for the details. The so-called \emph{enhanced virtual element space} in two dimensions is now defined as follows:
\begin{equation}
\mathbb{V}(F) := \left\{v\in \tilde{\mathbb{V}}(F) \middle|\ \int_F vp_1 = \int_F (\Pi_F^\nabla v)p_1,\ \forall\ p_1\in\mathbb{P}_1(F)\right\}.
\end{equation}

\noindent
For the three dimensional spaces, let now $E$ be a polyhedron.  The \emph{boundary space} on $\partial E$ and the preliminary virtual element space on $E$ are defined by
\begin{align*}
&\mathcal{B}(\partial E) = \{v \in \mathcal{C}^0(\partial E) | v_{|F} \in \mathbb{V}(F), \ \forall F \in \text{faces}(E)\};\\
&\tilde{\mathbb{V}}(E) := \Big\{v \in H^1(E) \Big|\ v_{|\partial E} \in\mathcal{B}(\partial E) \land \Delta v \in \mathbb{P}_1(E) \Big\},
\end{align*}
where $\mathbb{P}_1(E)$ is the space of linear polynomials on the polyhedron $E$. Let us consider the elliptic projection $\Pi^\nabla_E : \tilde{\mathbb{V}}(E) \rightarrow \mathbb{P}_1(E)$ defined by
\begin{equation}
\displaystyle\int_E \nabla (v - \Pi^\nabla_E v) \cdot \nabla p_1 = 0 \qquad \forall\ p_1 \in \mathbb{P}_1(E) \quad \land \quad
\displaystyle\int_{\partial E} (v - \Pi^\nabla_E v) = 0.
\end{equation}
Once again, the operator $\Pi_E^\nabla$ is computable, see \cite{ahmad2013equivalent} for the details. The enhanced virtual element space in three dimensions is now defined as follows:
\begin{equation}
\mathbb{V}(E) := \left\{v\in \tilde{\mathbb{V}}(E) \middle|\ \int_E vp_1 = \int_E (\Pi_E^\nabla v)p_1,\ \forall\ p_1\in\mathbb{P}_1(E)\right\}.
\end{equation}
The practical usability of the spaces $\mathbb{V}(F)$ and $\mathbb{V}(E)$ stem from the following result.
\begin{proposition}[Degrees of freedom (\cite{ahmad2013equivalent})]
\label{prop:dofs}
Let $n\in\mathbb{N}$. If $E$ is a polygon or a polyhedron with $n$ vertexes $\boldx_i$, $i=1,\dots,n$, then $\dim(\mathbb{V}(E)) = n$ and each function $v \in \mathbb{V}(E)$ is uniquely defined by the nodal values $v(\boldx_i)$, $i=1,\dots,n$. Hence, the nodal values constitute a set of degrees of freedom.
\end{proposition}

The following definition allows to correctly handle functions that are multiply defined on the junction between elements. 
\begin{definition}[Broken Sobolev norms]
Given two collections of functions\\ $\{u_E: E\rightarrow \mathbb{R} | E \in\mathcal{E}_h\}$ and $\{v_F: F\rightarrow\mathbb{R} | F \in\mathcal{F}_h\}$,  the broken Sobolev seminorms are defined as follows:
\begin{align*}
|u|_{s,\Omega,h} := \left(\sum_{E\in\mathcal{E}_h} |u_E|_{H^s(E)}^2\right)^{\frac{1}{2}}, \qquad
|v|_{s,\Gamma,h} := \left(\sum_{F\in\mathcal{F}_h} |v_F|_{H^s(F)}^2\right)^{\frac{1}{2}}, \qquad s=1,2.
\end{align*}
\end{definition}

The approximation properties of the spaces $\mathbb{V}(F)$ and $\mathbb{V}(E)$ are given by the following result.
\begin{proposition}[Projection error on polynomials (\cite{ahmad2013equivalent})]
For $s=1,2$, given two collections of functions $\{u_E \in H^s(E)| E \in\mathcal{E}_h\}$ and $\{v_F \in H^s(F) | F \in\mathcal{F}_h\}$,  there exist $u_\pi \in \prod_{E\in\mathcal{E}_h}\mathbb{P}_1(E)$ and $v_\pi \in \prod_{F\in\mathcal{F}_h}\mathbb{P}_1(F)$ such that
\begin{align}
\label{projection_error}
&\|u-u_\pi\|_{0, \Omega,h} + h |u-u_\pi|_{1,\Omega,h} \leq C h^s |u|_{s,\Omega,h};\\
\label{projection_error_surf}
&\|v-v_\pi\|_{0,\Gamma,h} + h |v-v_\pi|_{1,\Gamma,h} \leq C h^s |v|_{s,\Gamma,h},
\end{align}
where $C$ is a constant that depends only on $\gamma_1$.
\end{proposition}

The \emph{global virtual element spaces} in the bulk and on the surface are defined by matching of degrees of freedom across elements:
\begin{align}
&\mathbb{V}_{\Omega} := \{v \in H^1(\Omega_h)\  | \ v_{|E} \in \mathbb{V}(E), \ \forall\ E \in \mathcal{E}_h\};\\
&\mathbb{V}_{\Gamma} := \{v \in \mathcal{C}^0(\Gamma_h) \ | \ v_{|F} \in \mathbb{V}(F), \ \forall\ F \in \mathcal{F}_h\}.
\end{align}
In the global spaces $\mathbb{V}_\Gamma$ and $\mathbb{V}_\Omega$ we consider the Lagrange basis functions $\varphi_{i} \in \mathbb{V}_\Omega$ for $i=1,\dots,N$ and $\psi_{i'} \in \mathbb{V}_\Gamma$ for $i'=1,\dots,M$, where each $\varphi_i$ and each $\psi_{i'}$ are uniquely defined by $\varphi_i(\boldx_j) = \delta_{ij}$ for all $i,j = 1,\dots,N$ and $\psi_{i'}(\boldx_{j'}) = \delta_{i'j'}$ for all $i',j' = 1,\dots,M$, respectively,
with $\delta_{ij}$ being the Kronecker symbol. The sets $\{\varphi_i, \ i=1,\dots,N\}$ and $\{\psi_{i'}, \ i'=1,\dots,M\}$ are bases of $\mathbb{V}_\Omega$ and $\mathbb{V}_\Gamma$, respectively, thanks to Proposition \ref{prop:dofs}. It is easy to see that the bulk- and surface- Lagrange basis functions fulfil the following relation:
\begin{equation}
\label{basis_functions_relation}
\varphi_{i | \Gamma_h} = \psi_i, \qquad \forall\ i=1,\dots,M.
\end{equation}

\subsection{Approximation of bilinear forms}
In order to derive a spatially discrete formulation of the weak continuous problem  \eqref{elliptic_problem_weak_form} we need suitable approximate bilinear forms. We will follow \cite{ahmad2013equivalent, beirao2013basic}. In the remainder of this section, let $F$ and $E$ be elements of $\Gamma_h$ and $\Omega_h$, respectively. The \emph{stabilizing forms} $S_F: \mathbb{V}(F) \times \mathbb{V}(F) \rightarrow \mathbb{R}$ and $S_E: \mathbb{V}(E) \times \mathbb{V}(E) \rightarrow \mathbb{R}$ are defined by
\begin{align}
&S_F(v,w) := \sum_{P \in \text{ vertexes }(F)} v(P)w(P), \qquad \forall \ v,\, w\in \mathbb{V}(F);\\
&S_E(v,w) := \sum_{P \in \text{ vertexes }(E)} v(P)w(P), \qquad \forall \ v,\, w\in \mathbb{V}(E),
\end{align}
respectively. The $L^2$ projectors $\Pi^0_F: \mathbb{V}(F) \rightarrow\mathbb{P}_1(F)$ and $\Pi^0_E: \mathbb{V}(E) \rightarrow\mathbb{P}_1(E)$ are defined as follows: for $v \in \mathbb{V}(F)$ and $w \in \mathbb{V}(E)$:
\begin{align}
\label{L2_projector_surf}
&\int_F (v - \Pi^0_F v)p = 0, \qquad \forall\ \ p\in\mathbb{P}_1(F);\\
\label{L2_projector}
&\int_E (w - \Pi^0_E w)p = 0, \qquad \forall\ \ p\in\mathbb{P}_1(E),
\end{align}
respectively. As shown in \cite{ahmad2013equivalent}, $\Pi_F^0$ and $\Pi_E^0$ are computable because $\Pi_F^0 = \Pi_F^\nabla$ and $\Pi_E^0 = \Pi_E^\nabla$. Even if $\Pi_F^0$ and $\Pi_E^0$ are not new projectors, the presentation and the analysis of the method benefit from the usage of the equivalent definitions \eqref{L2_projector_surf}-\eqref{L2_projector}. Moreover, since $\Pi_F^0 = \Pi_F^\nabla$ and $\Pi_E^0 = \Pi_E^\nabla$, the boundedness property of projection operators in Hilbert spaces translates to
\begin{align}
\label{L2_projector_boundedness_surf}
&\|\Pi^0_F v\|_{L^2(F)} \leq \|v\|_{L^2(F)} \quad \text{and} \qquad |\Pi^0_F v|_{H^1(F)} \leq |v|_{H^1(F)};\\
\label{L2_projector_boundedness}
&\|\Pi^0_E w\|_{L^2(E)} \leq \|w\|_{L^2(E)} \quad \text{and} \qquad |\Pi^0_E w|_{H^1(E)} \leq |w|_{H^1(E)}.
\end{align}
We are now ready to introduce the approximate $L^2$ bilinear forms $m_F: \mathbb{V}(F) \times \mathbb{V}(F) \rightarrow\mathbb{R}$ and $m_E: \mathbb{V}(E) \times \mathbb{V}(E) \rightarrow\mathbb{R}$, defined as follows:
\begin{align}
&m_F(v,w) := \int_F (\Pi^0_F v)(\Pi^0_F w) + \area(F) S_F(v-\Pi^0_F v, w - \Pi^0_F w);\\
&m_E(v,w) := \int_E (\Pi^0_E v)(\Pi^0_E w) + \volume(E) S_E(v-\Pi^0_E v, w - \Pi^0_E w),
\end{align}
respectively. The approximate gradient-gradient bilinear forms $a_F: \mathbb{V}(F) \times  \mathbb{V}(F) \rightarrow\mathbb{R}$ and $a_E: \mathbb{V}(E) \times  \mathbb{V}(E) \rightarrow\mathbb{R}$ are defined by
\begin{align}
&a_F(v,w) := \int_F (\nabla \Pi^\nabla_F v) \cdot (\nabla \Pi^\nabla_F w) + \diam(F) S_F(v-\Pi^\nabla_F v, w - \Pi^\nabla_F w);\\
&a_E(v,w) := \int_E (\nabla \Pi^\nabla_E v) \cdot (\nabla \Pi^\nabla_E w) + \diam(E) S_E(v-\Pi^\nabla_E v, w - \Pi^\nabla_E w),
\end{align}
respectively. The definitions of $a_E$, $a_F$,  $m_E$ and $m_F$ imply the following result.
\begin{proposition}[Stability and consistency (\cite{beirao2013basic})]
\label{prop:stab_cons}
The bilinear forms $a_E$, $a_F$,  $m_E$ and $m_F$ are consistent, i.e.
\begin{align}
\label{consistency_surf}
&a_F(v,p) = \int_F \nabla v \cdot \nabla p; \quad m_F(v,p) = \int_F vp, \qquad \forall p \in \mathbb{P}_1(F);\\
\label{consistency}
&a_E(v,p) = \int_E \nabla v \cdot \nabla p; \quad m_E(v,p) = \int_E vp, \qquad \forall p \in \mathbb{P}_1(E).
\end{align}
The bilinear forms $a_E$, $a_F$,  $m_E$ and $m_F$ are stable, meaning that there exist two constants $0 < \alpha_* < \alpha^*$ depending on $\gamma_2$ such that, for all $v\in\mathbb{V}(F)$ and $w\in\mathbb{V}(E)$
\begin{align}
\label{stability_surf}
&\alpha_* \int_F \nabla v \cdot \nabla v \leq a_F(v,v) \leq \alpha^* \int_F \nabla v \cdot \nabla v; \quad  \alpha_* \int_F v^2 \leq m_F(v,v) \leq \alpha^* \int_F v^2;\\
\label{stability_bulk}
&\alpha_* \int_E \nabla w \cdot \nabla w \leq a_E(w,w) \leq \alpha^* \int_E \nabla w\cdot \nabla w; \quad \alpha_* \int_E w^2 \leq m_E(w,w) \leq \alpha^* \int_E w^2.
\end{align}
\end{proposition}
We observe from \eqref{stability_surf}-\eqref{stability_bulk} that the approximate bilinear forms $a_E$, $a_F$,  $m_E$ and $m_F$ do not converge to their exact counterparts,  see also \cite{beirao2013basic}. Nevertheless, we will show that the method retains optimal convergence thanks to the consistency properties \eqref{consistency_surf}-\eqref{consistency}.
The global bilinear forms $a_h^\Gamma, m_h^\Gamma: \mathbb{V}_\Gamma \times \mathbb{V}_\Gamma \rightarrow\mathbb{R}$, and $a_h^\Omega, m_h^\Omega: \mathbb{V}_\Omega \times \mathbb{V}_\Omega \rightarrow\mathbb{R}$ are defined elementwise:
\begin{align}
&a_h^\Gamma(v,w) := \sum_{F\in\mathcal{F}_h}  a_F(v_{|F}, w_{|F}); \quad m_h^\Gamma(v,w) := \sum_{F\in\mathcal{F}_h}  m_F(v_{|F}, w_{|F});\\
&a_h^\Omega(v,w) := \sum_{E\in\mathcal{E}_h}  a_E(v_{|E}, w_{|E}); \quad m_h^\Omega(v,w) := \sum_{E\in\mathcal{E}_h}  m_E(v_{|E}, w_{|E}).
\end{align}
From Proposition \ref{prop:stab_cons}, $m_h^\Gamma$ and $m_h^\Omega$ are positive definite, while $a_h^\Gamma$ and $a_h^\Omega$ are positive semi-definite.

\subsection{Approximation of the load terms}
The approximate bilinear forms $m_h^\Gamma$ and $m_h^\Omega$ presented in the previous section are not sufficient to discretise load terms like $\intgamma g\varphi$ and $\intomega f\varphi$, because $g$ and $f$ are not in the spaces $\mathbb{V}_\Gamma$ and $\mathbb{V}_\Omega$, respectively. 

\begin{definition}[Surface- and bulk- virtual Lagrange interpolants]
Given $f \in\mathcal{C}^0(E)$, $E\in\mathcal{E}_h$ and $g \in \mathcal{C}^0(F)$, $F \in \mathcal{F}_h$, the \emph{virtual Lagrange interpolants} $I_\Omega f$ of $f$ and $I_\Gamma g$ of $g$ are the unique $\mathbb{V}(E)$ and $\mathbb{V}(F)$ functions, respectively, such that $I_E f(\boldx) = f(\boldx)$ for all $\boldx \in \text{ vertexes}(E)$ and $I_F g(\boldx) = g(\boldx)$ for all $\boldx \in \text{ vertexes}(F)$, respectively.  Given two collections of functions $\{f_E \in \mathcal{C}^0(E) | E \in\mathcal{E}_h\}$ and $\{g_F \in \mathcal{C}^0(F) | F \in\mathcal{F}_h\}$, their \emph{global interpolants} are the collections of functions defined by $I_\Omega f = \{I_E f_E | E \in\mathcal{E}_h\}$ and $I_\Gamma g = \{I_F g_F | F \in\mathcal{F}_h\}$. 
\end{definition}

\begin{proposition}[Interpolation error (\cite{ahmad2013equivalent})]
\label{prop:interpolation_error}
Given two collections of functions $\{f_E \in H^2(E) | E \in\mathcal{E}_h\}$ and $\{g_F \in H^2(F) | F \in\mathcal{F}_h\}$,  it holds that
\begin{align}
\label{interpolation_error_bulk}
&\|f-I_\Omega(f)\|_{L^2(\Omega_h)} + h|f-I_\Omega(f)|_{1,\Omega,h} \leq Ch^2|f|_{2,\Omega,h};\\
\label{interpolation_error_surf}
&\|g-I_\Gamma(g)\|_{L^2(\Gamma_h)} + h|g-I_\Gamma(g)|_{1,\Gamma,h} \leq Ch^2|g|_{2,\Gamma,h},
\end{align}
respectively, where $C>0$ depends only on $\gamma_1$.
\end{proposition}

\subsection{The spatially discrete formulation}
The discrete counterpart of the elliptic problem \eqref{elliptic_problem_weak_form} is: find $(U,V) \in \mathbb{V}_\Omega \times \mathbb{V}_\Gamma$ such that
\begin{equation}
\label{elliptic_problem_BSVEM}
b_h((U,V); (\varphi,\psi)) = m_h^\Omega(I_\Omega(f), \varphi) + m_h^\Gamma(I_\Gamma(g), \psi), \quad \forall (\varphi,\psi) \in \mathbb{V}_\Omega \times \mathbb{V}_{\Gamma},
\end{equation}
where $b_h: (\mathbb{V}_\Omega \times \mathbb{V}_\Gamma)^2 \rightarrow\mathbb{R}$ is the discrete bilinear form defined by
\begin{equation}
\label{discrete_bilinear_form}
\begin{split}
b_h((U,V);(\varphi,\psi)) &:= \alpha\Big(a_h^\Omega(U,\varphi) + m_h^\Omega(U,\varphi)\Big) + \beta\Big(a_h^\Gamma(V,\psi) + m_h^\Gamma(V,\psi)\Big)\\
&+ m_h^\Gamma(\alpha U - \beta V, \alpha \varphi - \beta \psi).
\end{split}
\end{equation}
We express the spatially discrete solution $(U,V)$ in the Lagrange bases as follows:
\begin{align}
U(\boldx) = \sum_{i=1}^N \xi_i \varphi_i(\boldx), \qquad \boldx \in \Omega_h;\quad \text{and} \quad 
V(\boldx) = \sum_{k=1}^M \eta_k \psi_k(\boldx), \qquad \boldx \in \Gamma_h.
\end{align}
Hence, problem \eqref{elliptic_problem_BSVEM} is equivalent to: find $\boldxi := (\xi_i, \dots, \xi_N)^{T} \in \mathbb{R}^{N}$ and $\boldeta := (\eta_1,\dots, \eta_M)^{T} \in\mathbb{R}^{M}$ such that 
\begin{equation}
\label{elliptic_problem_BSVEM_basis}
\begin{cases}
\begin{aligned}
\displaystyle\sum_{i=1}^N\displaystyle \xi_i  \Big(a_h^\Omega(\varphi_i, \varphi_j) +  m_h^\Omega(\varphi_i, \varphi_j)\Big)  +  &\sum_{k=1}^M \Big( \alpha \xi_k m_h^\Gamma(\varphi_k, \varphi_l)  - \beta  \eta_k m_h^\Gamma(\psi_k, \varphi_l)\Big) \\
= &\sum_{i=1}^N f(\boldx_i) m_h^\Omega(\varphi_i, \varphi_j) \qquad \forall j=1,\dots,N;
\end{aligned}\\
\begin{aligned}
\displaystyle \sum_{k=1}^M \Big( \eta_k a_h^\Gamma(\psi_k, \psi_l) -\alpha  \xi_k m_h^\Gamma(\varphi_k, \psi_l) &+ (\beta+1)  \eta_k m_h^\Gamma(\psi_k, \psi_l) \Big)\\
 & \hspace*{-6mm}= \sum_{k=1}^M g(\boldx_k)m_h^\Gamma(\psi_k, \psi_l) \qquad \forall l=1,\dots,M.
\end{aligned}
\end{cases}
\hspace*{-5mm}
\end{equation}
We define the matrices $A_\Omega = (a_{i,j}^{\Omega}) \in\mathbb{R}^{N\times N}$,  $M_\Omega = (m_{i,j}^{\Omega}) \in\mathbb{R}^{N \times N}$, $A_\Gamma = (a_{k,l}^{\Gamma}) \in \mathbb{R}^{M \times M}$ and $M_{\Gamma} = (m_{k,l}^{\Gamma}) \in\mathbb{R}^{M \times M}$ as follows:
\begin{align}
\label{definition_matrices_bulk}
&a_{i,j}^{\Omega} := a_h^\Omega(\varphi_i, \varphi_j),  \quad \text{and} \quad m_{i,j}^{\Omega} := m_h^\Omega(\varphi_i, \varphi_j), \qquad i,j=1,\dots,N;\\
\label{definition_matrices_surf}
&a_{k,l}^{\Gamma} := a_h^\Gamma(\psi_k, \psi_l),  \quad \text{and} \quad m_{k,l}^{\Gamma} := m_h^\Gamma(\psi_k, \psi_l), \qquad k,l=1,\dots,M.
\end{align}
By using \eqref{basis_functions_relation} and defining $\boldf := (f(\boldx_1), \dots, f(\boldx_N))^T \in\mathbb{R}^{N}$ and $\boldg := (g(\boldx_1), \dots, g(\boldx_M))^T \in \mathbb{R}^{M}$ we can rewrite the discrete formulation \eqref{elliptic_problem_BSVEM_basis} as a block $(N+M)\times(N+M)$ linear algebraic system:
\begin{equation}
\label{elliptic_problem_BSVEM_linear_system}
\begin{cases}
A_\Omega \boldxi + M_\Omega \boldxi + \alpha RM_\Gamma R^T\boldxi - \beta R M_\Gamma \boldeta = M_\Omega \boldf;\\
A_\Gamma \boldeta - \alpha M_\Gamma R^T\boldxi + (\beta + 1) M_\Gamma \boldeta = M_\Gamma \boldg.
\end{cases}
\end{equation}
In compact form, the linear system \eqref{elliptic_problem_BSVEM_linear_system} reads
\begin{equation}
 \label{elliptic_problem_BSVEM_linear_system_compact}
\left[\begin{matrix}
A_\Omega + M_\Omega + \alpha RM_\Gamma R^T     & - \beta R M_\Gamma\\
- \alpha M_\Gamma R^T                                            & A_\Gamma + (\beta + 1) M_\Gamma
\end{matrix}\right]
\left[\begin{matrix}
\boldxi\\
\boldeta
\end{matrix}\right]
=
\left[\begin{matrix}
M_\Omega\boldf\\
M_\Gamma\boldg
\end{matrix}\right].
\end{equation}
It is possible to show that the coefficient matrix of \eqref{elliptic_problem_BSVEM_linear_system_compact} is sparse and unstructured.

\section{Convergence analysis}
\label{sec:convergence_analysis}

To derive error estimates for the  discrete solution we need a bulk-surface Ritz projection tailored for the variational problem \eqref{elliptic_problem_weak_form}.
\begin{definition}[Bulk-Surface Ritz projection]
The \emph{bulk-surface Ritz projection} of a pair $(u,v)\in H^1(\Gamma)\times H^1(\Omega)$ is the unique pair $(\mathcal{R}u, \mathcal{R}v) \in \mathbb{V}_\Omega \times \mathbb{V}_\Gamma$ such that
\begin{equation}
\label{bulk_surface_ritz}
b_h((U,V);(\varphi,\psi)) = b((u,v);(\varphi^\ell,\psi^\ell)), \qquad \forall (\varphi,\psi) \in \mathbb{V}_\Omega \times \mathbb{V}_\Gamma.
\end{equation}
\end{definition}

The bulk-surface Ritz projection is well-defined since $b_h$ is coercive.

\begin{theorem}[$H^1(\Omega) \times H^1(\Gamma)$ a priori error bound for the bulk-surface Ritz projection]
\label{thm:ritz_proj_H1}
For any $(u,v)\in H^{2+3/4}(\Omega)\times H^2(\Gamma)$ it holds that
\begin{equation}
\label{ritz_estimate_H1}
\begin{split}
\|(u,v) &- (\mathcal{R}u, \mathcal{R}v)^\ell\|_{H^1(\Omega)\times H^1(\Gamma)} \leq Ch \|(u,v)\|_{H^{2+3/4}(\Omega)\times H^2(\Gamma)},
\end{split}
\end{equation}
where the additional index $3/4$ appears only in the simultaneous presence of curved boundaries and non-tetrahedral exterior elements.

\begin{proof}
We set $e_h = (e_h^\Omega, e_h^\Gamma) := (\mathcal{R}u - u^{-\ell}, \mathcal{R}v - v^{-\ell})$. From \eqref{consistency}, \eqref{stability_surf}, \eqref{stability_bulk} and \eqref{bulk_surface_ritz} we have
\begin{align}
\label{ritz_h1_error_equation}
\alpha_*\min(\alpha,\beta)&\|e_h\|_{H^1(\Omega) \times H^1(\Gamma)} \leq b_h(e_h,e_h) = \alpha \Big( \underset{T_1}{\underbrace{a_h^\Omega(e_h^\Omega, e_h^\Omega)}} +  \underset{T_2}{\underbrace{m_h^\Omega(e_h^\Omega, e_h^\Omega)}} \Big) \notag\\
+ &\beta \Big( \underset{T_3}{\underbrace{ a_h^\Gamma(e_h^\Gamma, e_h^\Gamma)}} + \underset{T_4}{\underbrace{m_h^\Gamma(e_h^\Gamma, e_h^\Gamma)}} \Big) + \underset{T_5}{\underbrace{m_h^\Gamma(\alpha e_h^\Omega - \beta e_h^\Gamma,  \alpha e_h^\Omega - \beta e_h^\Gamma)}}.
\end{align}
We estimate $T_1$ using \eqref{lifting_error_bulk_nabla}, \eqref{error_between_lift_and_extension}, \eqref{projection_error}, \eqref{consistency}, \eqref{bulk_surface_ritz} and the continuity of $a_h^\Omega$:
\begin{align}
\label{ritz_projection_energy_T1}
T_1 &= a_h^\Omega(\mathcal{R}u, e_h) - a_h^\Omega(u^{-\ell}, e_h)\\
= \intomega \nabla u\cdot\nabla e_h^{\ell} &+ a_h^\Omega(\tilde{u} - u^{-\ell}, e_h) + a_h^\Omega( \tilde{u}_\pi - \tilde{u}, e_h) - a_h^\Omega(\tilde{u}_\pi,e_h) \notag\\
 = \intomega \nabla u\cdot\nabla e_h^{\ell} &- \intomegah \nabla \tilde{u}_\pi\cdot\nabla e_h + a_h^\Omega(\tilde{u} - u^{-\ell}, e_h) + a_h^\Omega( \tilde{u}_\pi - \tilde{u}, e_h) \notag\\
= \intomega \nabla u \cdot \nabla e_h^{\ell} &- \intomegah \nabla u^{-\ell}\cdot \nabla e_h + \intomegah \nabla (u^{-\ell} - \tilde{u})\cdot \nabla e_h  + \intomegah \nabla (\tilde{u} - \tilde{u}_\pi)\cdot \nabla e_h  \notag\\
 & + a_h^\Omega(\tilde{u} - u^{-\ell}, e_h) + a_h^\Omega( \tilde{u}_\pi - \tilde{u}, e_h) \notag\\
\leq  C\Big(h\|u\|_{H^2(\Omega)} &+ h^{\frac{3}{2}}\|u\|_{H^2(\Omega)} + Ch\|u\|_{H^{2+3/4}(\Omega)}\Big) |e_h|_{H^1(\Omega_h)} \notag.
\end{align}
We estimate $T_2$ in the same way by using \eqref{lifting_error_bulk_mass}, \eqref{error_between_lift_and_extension}, \eqref{consistency},\eqref{projection_error} and the continuity of $m_h^\Omega$:
\begin{align}
\label{ritz_projection_energy_T2}
T_2 \leq Ch^2\left(\|u\|_{H^2(\Omega)} + \|u\|_{H^{2+1/4}(\Omega)}\right) \|e_h\|_{L^2(\Omega_h)}.
\end{align}
We estimate $T_3$ by reasoning as for $T_1$, but this time there is no need for the Sobolev extension because, as opposed to the $H^2(\Omega)$ norm,  the $H^2(\Gamma)$ norm is preserved under lifting thanks to \eqref{equivalence_gamma_h2}. This implies that $v^{-\ell}$ is $H^2$ on each face of $\Gamma_h$ and thus fulfils the optimal error estimate for the projection \eqref{projection_error_surf}.  Hence, by using \eqref{lifting_error_surf_nabla}, \eqref{consistency_surf} and \eqref{projection_error_surf}, the estimate for $T_3$ reads as follows:
\begin{align}
T_3 = a_h^\Gamma(\mathcal{R}v, e_h^\Gamma) - a_h^\Gamma(v^{-\ell}, e_h^\Gamma) =  \intgamma \nablagamma v\cdot\nablagamma e_h^{\Gamma, \ell} + a_h^\Gamma(v^{-\ell}_\pi - v^{-\ell}, e_h) - a_h^\Gamma(v^{-\ell}_\pi,e_h^\Gamma) \notag\\
 = \intgamma \nablagamma v\cdot\nablagamma e_h^{\Gamma, \ell} - \intgammah \nablagammah v^{-\ell}_\pi\cdot\nablagammah e_h^{\Gamma} + a_h^\Gamma(v^{-\ell}_\pi - v^{-\ell}, e_h^\Gamma) \notag\\
= \intgamma \nablagamma v\cdot\nablagamma e_h^{\Gamma, \ell} - \intgammah \nablagammah v^{-\ell}\cdot\nablagammah e_h^{\Gamma} + \intgammah \nablagammah (v^{-\ell} - v^{-\ell}_\pi)\cdot\nablagammah e_h^{\Gamma} \notag\\
\label{ritz_projection_energy_T3}
+ a_h^\Gamma(v^{-\ell}_\pi - v^{-\ell}, e_h^\Gamma) \leq  Ch\|v\|_{H^2(\Gamma)} |e_h|_{H^1(\Gamma_h)}. 
\end{align}
We estimate $T_4$ in the same way as $T_3$, by using \eqref{lifting_error_surf_mass} instead of \eqref{lifting_error_surf_nabla} and choosing $s=1$ instead of $s=2$ in \eqref{projection_error}:
\begin{align}
\label{ritz_projection_energy_T4}
T_4 \leq & Ch\|v\|_{H^1(\Gamma)} \|e_h\|_{L^2(\Gamma_h)}.
\end{align}
We estimate $T_5$ exactly as $T_4$ and then we apply the inverse trace inequality \eqref{inverse_trace_inequality}:
\begin{align}
\label{ritz_projection_energy_T5}
T_5 &\leq Ch\Big(\|v\|_{H^1(\Gamma)} + \|\Tr u\|_{H^1(\Gamma)}\Big) \Big( \|e_h^\Gamma\|_{L^2(\Gamma_h)} + \|\Tr e_h^\Omega\|_{L^2(\Gamma_h)}\Big) \notag\\
&\leq Ch\Big(\|v\|_{H^1(\Gamma)} + \|u\|_{H^2(\Gamma)}\Big)  \Big( \|e_h^\Gamma\|_{L^2(\Gamma_h)} + \|e_h^\Omega\|_{H^1(\Gamma_h)}\Big).
\end{align}
By substituting \eqref{ritz_projection_energy_T1}-\eqref{ritz_projection_energy_T5} into \eqref{ritz_h1_error_equation} and applying a Young inequality argument, we get the desired estimate \eqref{ritz_estimate_H1}.
In \eqref{ritz_projection_energy_T1}-\eqref{ritz_projection_energy_T2} notice that, in the absence of curvature or non-tetrahedral exterior elements, $u^{-\ell}_{|E} \in H^2(E)$ for all elements $E\in\mathcal{E}_h$, see Remark \ref{rmk:lifting_regularity}. Then the Sobolev extension $\tilde{u}$ is not needed and the terms in $H^{2+3/4}(\Omega)$ and $H^{2+1/4}(\Omega)$ do not appear. This completes the proof.
\end{proof}
\end{theorem}

\begin{theorem}[$L^2(\Omega) \times L^2(\Gamma)$ error bound for the bulk-surface Ritz projection]
\label{thm:ritz_proj_L2}
Let $\Omega$ have a $\mathcal{C}^3$ boundary. Then, for any $(u,v)\in H^{2+3/4}(\Omega)\times H^2(\Gamma)$ and for $h$ sufficiently small, it holds that
\begin{equation}
\label{bulk_ritz_estimate_bulk_L2}
\begin{split}
\|(u,v) &- (\mathcal{R}u, \mathcal{R}v)^\ell\|_{L^2(\Omega)\times L^2(\Gamma)} \leq Ch^2 \|(u,v)\|_{H^{2+3/4}(\Omega)\times H^2(\Gamma)},
\end{split}
\end{equation}
with $C$ depending on $\Omega$, $\gamma_1$ and $\gamma_2$. In \eqref{bulk_ritz_estimate_bulk_L2}, the additional exponent $3/4$ arises only in the simultaneous presence of curved boundaries and non-tetrahedral exterior elements. 

\begin{proof}
We will use an adapted Aubin-Nitsche duality method. Consider the dual problem: find $(\eta, \theta) \in H^1(\Omega) \times H^1(\Gamma)$ such that
\begin{equation}
\label{aubin_nitsche_problem}
b((\eta,\theta); (\varphi, \psi)) = \intomega (u-(\mathcal{R}u)^\ell) \varphi + \intgamma (v-(\mathcal{R}v)^\ell)\psi, 
\end{equation}
for all $(\varphi,\psi) \in H^1(\Omega) \times H^1(\Gamma)$. Since $u-(\mathcal{R}u)^\ell \in H^1(\Omega)$,  thanks to \eqref{exact_solution_bound_H^3}, the variational problem \eqref{aubin_nitsche_problem} has a unique solution $(\eta,\theta) \in H^3(\Omega) \times H^2(\Gamma)$ that fulfils
\begin{align}
\label{elliptic_regularity_2}
&\|(\eta, \theta)\|_{H^3(\Omega) \times H^2(\Gamma)} \leq C\|(u,v)-(\mathcal{R}u, \mathcal{R}v)^\ell)\|_{H^1(\Omega)\times H^1(\Gamma)}.
\end{align}
By combining \eqref{ritz_estimate_H1} and \eqref{elliptic_regularity_2} we have that
\begin{equation}
\label{elliptic_regularity_4}
\|(\eta, \theta)\|_{H^3(\Omega) \times H^2(\Gamma)} \leq Ch\|(u,v)\|_{H^2(\Omega) \times H^2(\Gamma)} + Ch\|u\|_{H^{2+3/4}(\Omega)}.
\end{equation}
We can choose $(\varphi,\psi) = (e_h^\Omega, e_h^\Gamma) = (u,v) - (\mathcal{R}u, \mathcal{R}v)^\ell$ in \eqref{aubin_nitsche_problem} and we get
\begin{equation}
\label{L2_bound_ritz_1}
\|(u,v) - (\mathcal{R} u, \mathcal{R}v)^\ell \|_{L^2(\Omega) \times L^2(\Gamma)}^2 = \intomega \varphi^2 + \intgamma \psi^2 = b((\eta,\theta); (u-\mathcal{R}u^\ell,  v - \mathcal{R}v^\ell)).
\end{equation}
The right hand side of \eqref{L2_bound_ritz_1} can be split into five terms, say $\bar{T}_1, \dots, \bar{T}_5$ as in \eqref{ritz_h1_error_equation}.  We explicitly show the estimation of the first of such terms -the most involved. The treatment of the other terms is similar. By using \eqref{lifting_error_bulk_nabla} and \eqref{bulk_surface_ritz} we have
\begin{equation}
\label{L2_bound_ritz_2}
\begin{split}
&\bar{T}_1 := a_h^\Omega(\eta, e_h^\Omega) = \intomega \nabla \eta\cdot \nabla \left(u - (\mathcal{R}u)^\ell\right)\\
= &\intomega \nabla \left(u - (\mathcal{R}u)^\ell\right)\cdot \nabla\left(\eta-I_\Omega(\tilde{\eta})^\ell\right) - \intomega \nabla (\mathcal{R}u)^\ell \cdot \nabla I_\Omega(\tilde{\eta})^\ell + a_h(\mathcal{R}u, I_\Omega(\tilde{\eta}))\\
\leq &|u - (\mathcal{R}u)^\ell|_{H^1(\Omega)}|\eta-I_\Omega(\tilde{\eta})^\ell|_{H^1(\Omega)} - \intomega \nabla (\mathcal{R}u)^\ell \cdot \nabla I_\Omega(\tilde{\eta})^\ell + a_h(\mathcal{R}u, I_\Omega(\tilde{\eta}))\\
\leq &|u - (\mathcal{R}u)^\ell|_{H^1(\Omega)}|\eta-I_\Omega(\tilde{\eta})^\ell|_{H^1(\Omega)} + Ch|(\mathcal{R}u)^\ell|_{H^1(\Omega_B^\ell)}|I_\Omega(\tilde{\eta})|_{H^1(\Omega_B^\ell)}\\
-&\intomegah \nabla \mathcal{R}u \cdot \nabla I_\Omega(\tilde{\eta}) + a_h(\mathcal{R}u, I_\Omega(\tilde{\eta}))\\
\leq &C\left(|u - (\mathcal{R}u)^\ell|_{H^1(\Omega)}  + h^\frac{1}{2}|u|_{H^1(\Omega_B^\ell)}\right)\left(|\eta-I_\Omega(\tilde{\eta})^\ell|_{H^1(\Omega)} + h^\frac{1}{2}|\eta|_{H^1(\Omega_B^\ell)}\right)\\
-&\intomegah \nabla \mathcal{R}u\cdot \nabla I_\Omega(\tilde{\eta}) + a_h^\Omega(\mathcal{R}u, I_\Omega(\tilde{\eta})),\\
\end{split}
\end{equation}
where we have used $h < h_0$ in the last inequality. We are left to estimate the right-hand-side of \eqref{L2_bound_ritz_2} piecewise. First, from \eqref{narrow_band_inequality} and \eqref{ritz_estimate_H1} we have that
\begin{equation}
\label{L2_bound_ritz_3}
\begin{split}
|u - (\mathcal{R}u)^\ell|_{H^1(\Omega)}  + h^\frac{1}{2}|u|_{H^1(\Omega_B^\ell)} \leq &Ch\|(u,v)\|_{H^2(\Omega) \times H^2(\Gamma)} + Ch\|u\|_{H^{2+3/4}(\Omega)}.
\end{split}
\end{equation}
Moreover, from \eqref{narrow_band_inequality}, \eqref{sobolev_extension},  \eqref{error_between_lift_and_extension}, \eqref{interpolation_error_bulk}, \eqref{elliptic_regularity_2} and \eqref{elliptic_regularity_4} we have that
\begin{equation}
\label{L2_bound_ritz_4}
\begin{split}
&|\eta-I_\Gamma(\tilde{\eta})^\ell|_{H^1(\Omega)} + h^\frac{1}{2}|\eta|_{H^1(\Omega_B^\ell)} \leq C|\eta^{-\ell}-I_\Gamma(\tilde{\eta})|_{H^1(\Omega_h)} + Ch\|\eta\|_{H^2(\Omega)}\\
\leq &C|\eta^{-\ell}-\tilde{\eta}|_{H^1(\Omega_h)} + C|\tilde{\eta}-I_\Gamma(\tilde{\eta})|_{H^1(\Omega_h)} + Ch\|\eta\|_{H^2(\Omega)}\\
\leq &Ch^{2}\|\eta\|_{H^3(\Omega)} + Ch\|\tilde{\eta}\|_{H^2(\Omega_h)} + Ch\|\eta\|_{H^2(\Omega)} \leq Ch\|\eta\|_{H^2(\Omega)} + Ch^{2}\|\eta\|_{H^3(\Omega)}\\
\leq &Ch\|u-(\mathcal{R}u)^\ell\|_{L^2(\Omega)} + Ch^3\|(u,v)\|_{H^2(\Omega) \times H^2(\Gamma)} + Ch^3\|u\|_{H^{2+3/4}(\Omega)}.
\end{split}
\end{equation}
Finally, we estimate the last two terms in \eqref{L2_bound_ritz_2} by adapting the approach used in \cite[Lemma 3.1]{vacca2015virtual}: from \eqref{error_between_lift_and_extension}, \eqref{projection_error}, \eqref{consistency}, \eqref{interpolation_error_bulk} and \eqref{elliptic_regularity_2} we have
\begin{equation}
\label{L2_bound_ritz_5}
\begin{split}
&a_h^\Omega(\mathcal{R}u, I_\Omega(\tilde{\eta})) - \intomegah \nabla \mathcal{R}u \cdot \nabla I_\Omega(\tilde{\eta})\\
= &\int_{\Omega_h} \nabla (\mathcal{R}u - \tilde{u}_\pi)\cdot\nabla(I_{\Omega}(\tilde{\eta}) - \tilde{\eta}_\pi) - a_h^\Omega(\mathcal{R}u - \tilde{u}_\pi, I_{\Omega}(\tilde{\eta}) - \tilde{\eta}_\pi)\\
\leq & |\mathcal{R}u - \tilde{u}_\pi|_{1,\Omega,h}|I_\Gamma(\tilde{\eta}) - \tilde{\eta}_\pi|_{1,\Omega,h} \leq Ch\|(u,v)\|_{H^{2+3/4}(\Omega) \times H^2(\Gamma)}Ch\|\eta\|_{H^2(\Omega)}\\
= &Ch^2\|(u,v)\|_{H^{2+3/4}(\Omega) \times H^2(\Gamma)}\|u-(\mathcal{R}u)^\ell \|_{L^2(\Omega)},
\end{split}
\end{equation}
By combining \eqref{L2_bound_ritz_2}-\eqref{L2_bound_ritz_5} we get
\begin{equation}
\label{bar_T_1_final_estimate}
\begin{split}
\bar{T}_1  \leq Ch^2 \|(u,v)\|_{H^{2+3/4}(\Omega) \times H^2(\Gamma)} \|e_h^\Omega\|_{L^2(\Omega)} + Ch^4 \|(u,v)\|_{H^{2+3/4}(\Omega) \times H^2(\Gamma)}^2.
\end{split}
\end{equation}
By estimating all remaining terms $\bar{T}_2\dots, \bar{T}_5$ as $\bar{T}_1$ in  \eqref{bar_T_1_final_estimate} and substituting into \eqref{L2_bound_ritz_1} we get
\begin{align}
\|(e_h^\Omega, e_h^\Gamma)\|_{L^2(\Omega)\times L^2(\Gamma)}^2  \leq &Ch^2 \|(u,v)\|_{H^{2+3/4}(\Omega) \times H^2(\Gamma)} \|(e_h^\Omega, e_h^\Gamma)\|_{L^2(\Omega)\times L^2(\Gamma)}\\
+ &Ch^4 \|(u,v)\|_{H^{2+3/4}(\Omega) \times H^2(\Gamma)}^2, \notag
\end{align}
where the additional index $3/4$ appears only in the simultaneous presence of curved boundaries and non-tetrahedral exterior elements, which proves \eqref{bulk_ritz_estimate_bulk_L2}.
\end{proof}
\end{theorem}

\begin{theorem}[$L^2(\Omega) \times L^2(\Gamma)$ error bound for the BSVEM]
\label{thm:BSVEM_L2_convergence}
Let $\Omega$ have a $\mathcal{C}^3$ boundary. Then,  if $(f,g) \in H^{2+1/4}(\Omega)\times H^2(\Gamma)$, the numerical solution $(U,V)$ fulfils
\begin{equation}
\label{BSVEM_error_estimate}
\begin{split}
\|(u,v) &- (U, V)^\ell\|_{L^2(\Omega)\times L^2(\Gamma)} \leq Ch^2 \|(f,g)\|_{H^{2+1/4}(\Omega)\times H^2(\Gamma)},
\end{split}
\end{equation}
with $C$ depending on $\Omega$, $\gamma_1$ and $\gamma_2$. In \eqref{BSVEM_error_estimate}, the additional index $1/4$ arises only in the simultaneous presence of curved boundaries and non-tetrahedral exterior elements.

\begin{proof}
The proof relies on a standard error equation technique.  The difference $\|(u,v) - (\mathcal{R}u, \mathcal{R}v)^\ell\|_{L^2(\Omega)\times L^2(\Gamma)}$ is estimated via \eqref{bulk_ritz_estimate_bulk_L2}, while the error equation for $\|(\mathcal{R}u, \mathcal{R}v)^\ell - (U,V)^\ell\|_{L^2(\Omega)\times L^2(\Gamma)}$,  obtained by subtracting the discrete problem \eqref{bulk_surface_ritz} from the weak continuous problem \eqref{elliptic_problem_weak_form},
is estimated via Lemma \ref{lmm:geometric_errors_bilinear forms},  Lemma \ref{lmm:geometric_error_of_sobolev_extension}, and Proposition  \ref{prop:interpolation_error}.
\end{proof}
\end{theorem}

\begin{remark}[Optimal convergence for bulk-only PDEs]
\label{rmk:bulk-only-case}
By considering the limit case $\beta = 0$ in the model problem \eqref{elliptic_problem}, the bulk equation becomes completely decoupled from the surface equation.  Specifically, the first equation in \eqref{elliptic_problem} becomes a linear elliptic equation in the 3D domain $\Omega$, endowed with non-zero Neumann boundary conditions.  Correspondingly, the BSVEM reduces to the known lowest-order VEM for elliptic problems in 3D (see \cite{beirao2013basic}).
Then, by setting $\beta=0$ throughout the present section devoted to convergence analysis, we obtain that the lowest-order VEM for 3D elliptic bulk problems retains optimal convergence in the presence of curved boundaries and non-zero boundary conditions. 
As mentioned in the Introduction, this result was not fully addressed in the literature, in this work we provide a rigorous justification.  It must be noted that previous works addressed the issue through the introduction of curved boundaries (see \cite{beirao2019curved,  dassi2021bend}) or the introduction of algebraic corrections in the method that account for surface curvature (see \cite{bertoluzza2019}). Here, we show for the first time that the plain 3D VEM of lowest order possesses optimal convergence even in the presence of curved boundaries.
\end{remark}

\section{Benefits of polyhedral meshes for BSPDEs}
\label{sec:mesh_advantage}
\begin{figure}[ht!]
\begin{center}
\subfigure[Step 1. The bulk $\Omega$, enclosed by the red boundary $\Gamma$, is bounded by the green cube $Q$, which is subdivided with a Cartesian grid.]{\label{fig:mesh_advantage_1}\includegraphics[scale=.4]{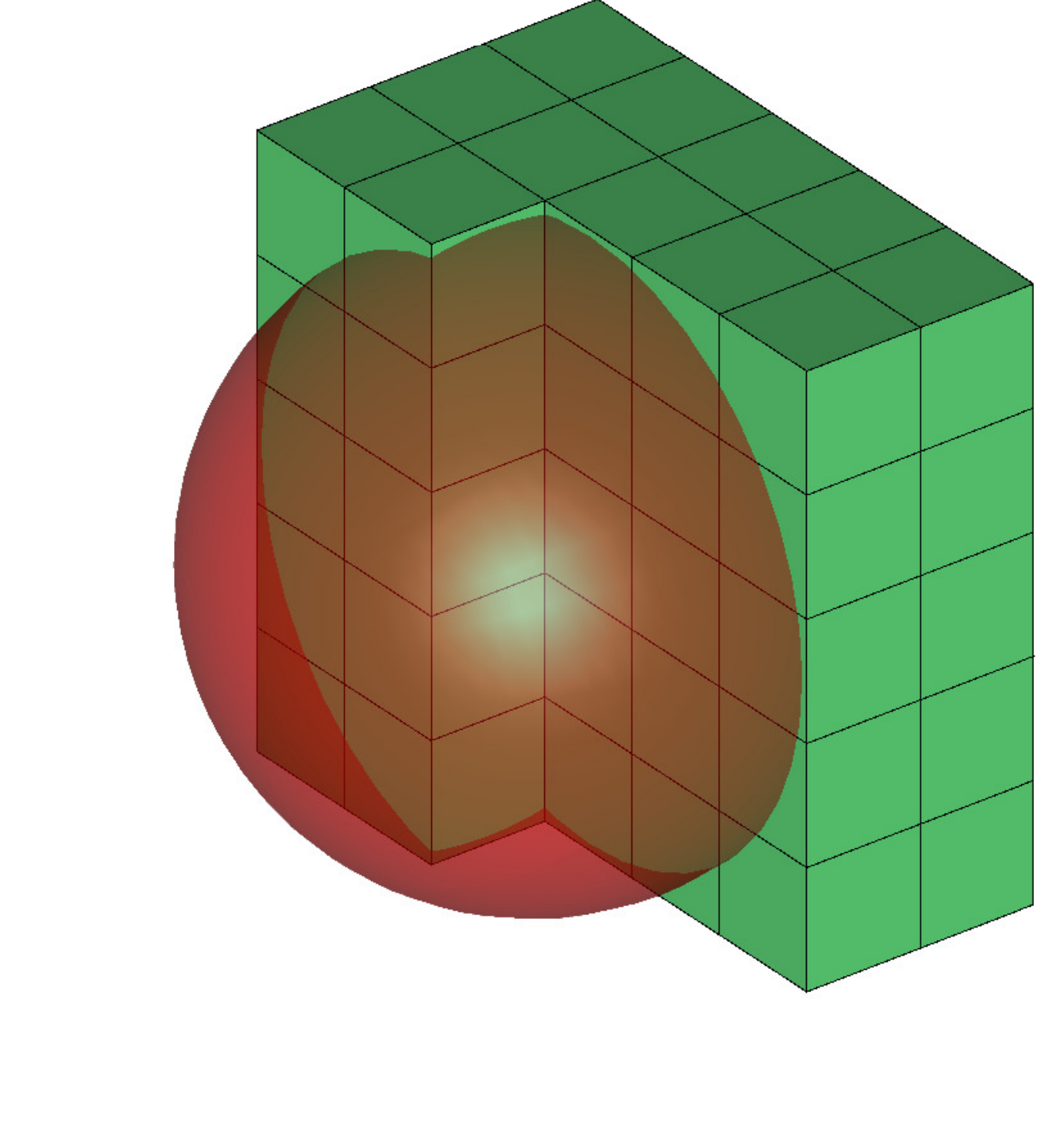}}
\subfigure[Step 2. The mesh elements that are not entirely inside $\Gamma$ are discarded.]{\label{fig:mesh_advantage_2}\includegraphics[scale=.4]{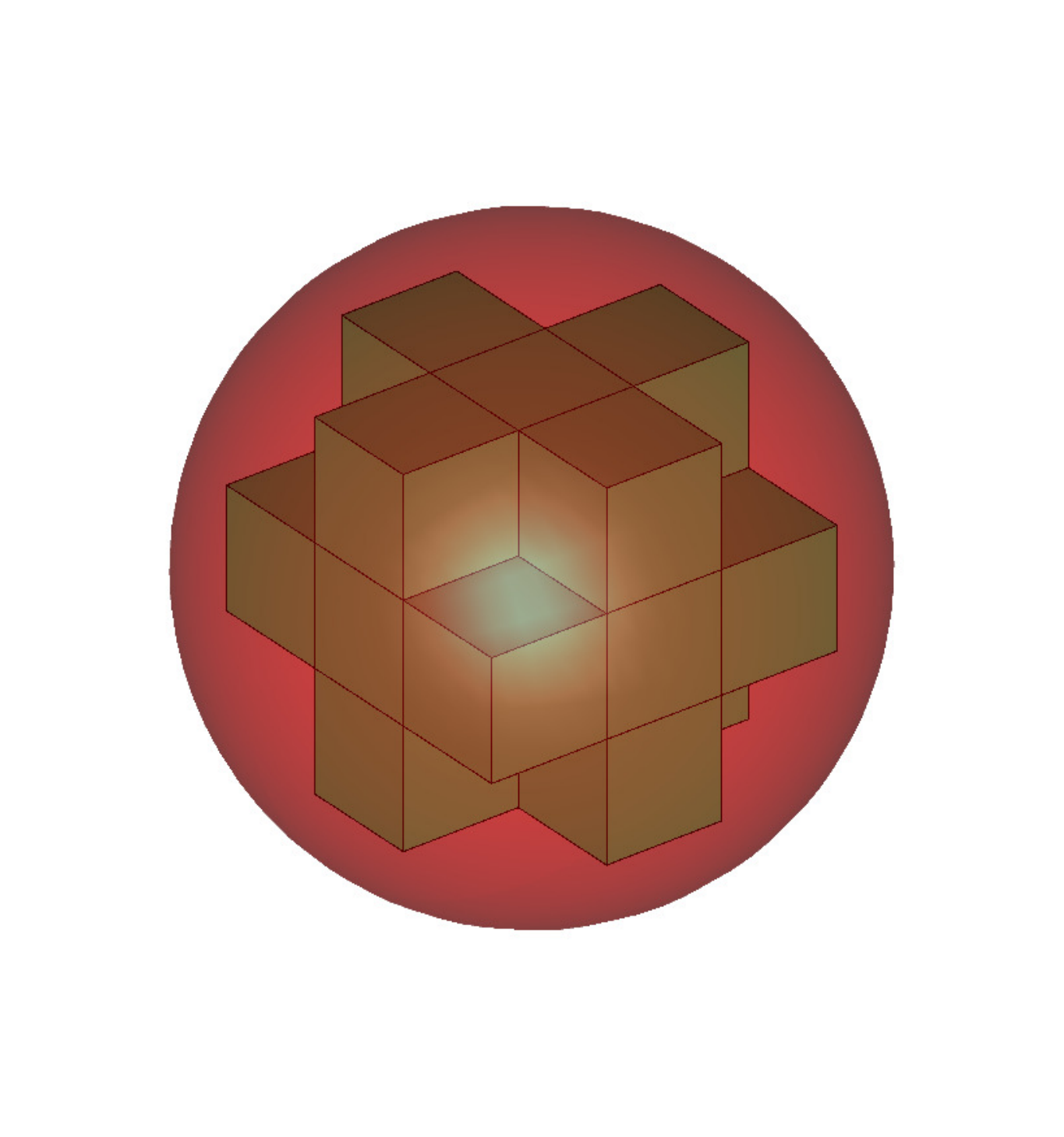}}
\subfigure[Step 3. The outermost cubic elements are \emph{extruded}, thereby producing the red band of polyhedral elements.]{\label{fig:mesh_advantage_3}\includegraphics[scale=.4]{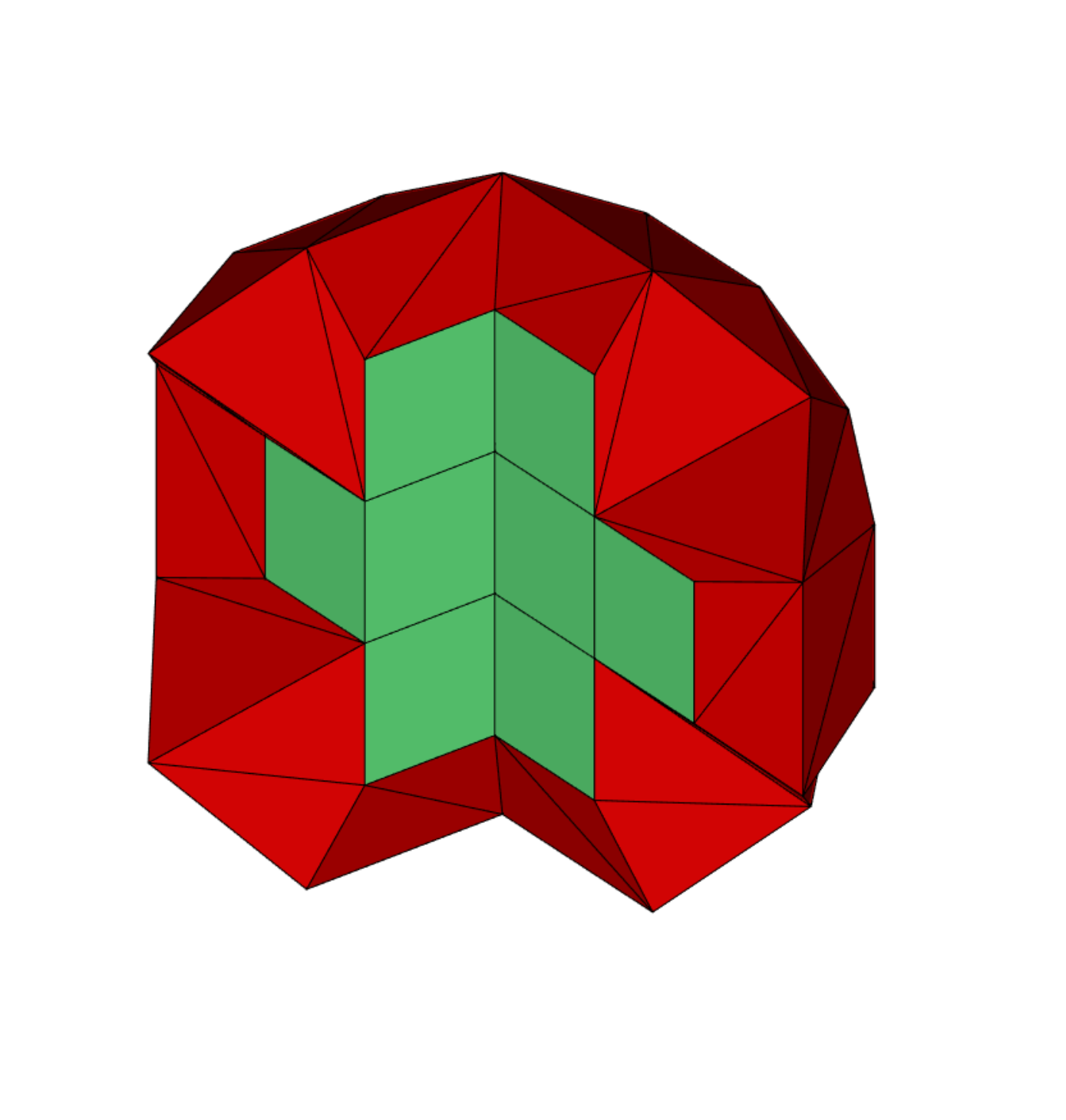}}
\subfigure[One of the polyhedral elements in Fig.  \ref{fig:mesh_advantage_3} generated through the extrusion process.]{\label{fig:mesh_advantage_4}\includegraphics[scale=.4]{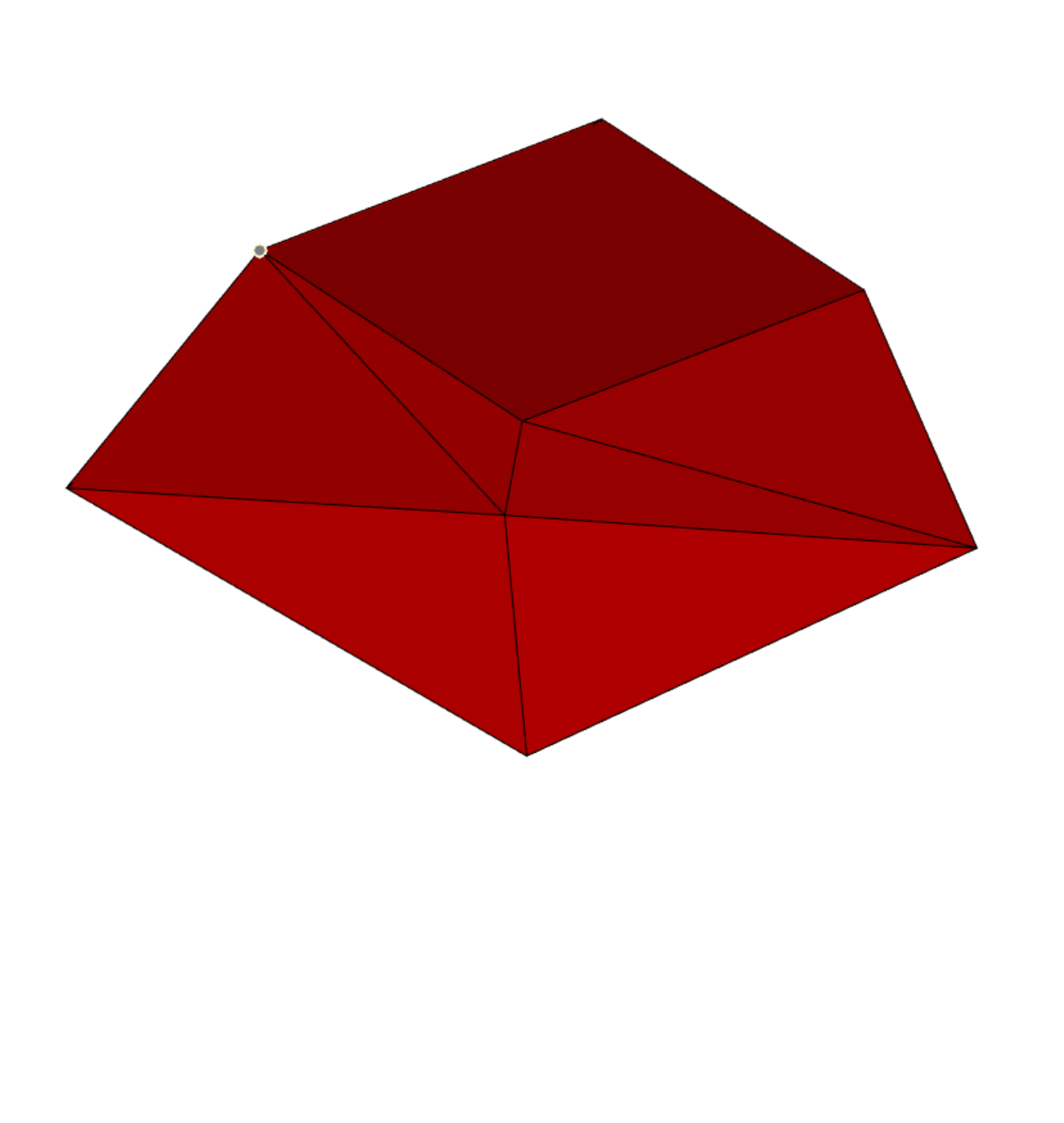}}
\end{center}
\caption{Generation of a polyhedral bulk-surface mesh that allows for optimised matrix assembly. }
\label{fig:mesh_advantage}
\end{figure}

If a domain $\Omega \subset \mathbb{R}^3$ has a $\mathcal{C}^1$ boundary $\Gamma$, we can construct a polyhedral mesh designed for fast matrix assembly, by proceeding as follows.  Enclose the bulk $\Omega$ in a cube $Q$. We discretise $Q$ with a Cartesian grid made up of cubic mesh elements and assume that at least one of such cubes is fully contained in $\Omega$ (see Fig.  \ref{fig:mesh_advantage_1}). Then we discard the elements that are not fully inside $Q$ (see Fig.  \ref{fig:mesh_advantage_2}), thereby producing an incomplete cubic mesh.  Finally,  we extrude the outermost (square) faces of the incomplete cubic mesh thus producing a discrete narrow band $\Omega_B$ of irregular polyhedral elements (highlighted in red in Fig. \ref{fig:mesh_advantage_3}). The resulting mesh $\Omega_h$ has the important property that it is made up of equal cubic elements, except for the exterior elements, as we can see in Fig. \ref{fig:mesh_advantage_3}.  This property allows for fast matrix assembly.  In fact if $h$ is the meshsize of $\Omega_h$, then the number on non-cubic elements of $\Omega_h$ is only $\mathcal{O}(h^{-2})$ out of $\mathcal{O}(h^{-3})$ overall elements, see \cite{frittelli2021bulk} for a discussion of the 2D case. This implies that, when assembling the mass- and stiffness- matrices $M_\Omega$ and $A_\Omega$ defined in \eqref{definition_matrices_bulk},  only $\mathcal{O}(h^{-2})$ element-wise local matrices must be actually computed, since the local matrices for a cubic element are known in closed form, see \cite{beirao2014hitchhiker}.

Matrix assembly optimization can be also achieved with different methods, such as cut FEM (\cite{burman2016cut}) or trace FEM (\cite{gross2015trace}). However, in these works, the authors adopt a level set representation of the boundary $\Gamma$, which we do not need in this study, as we exploit the usage of arbitrary polygons to approximate $\Gamma$.  Moreover, the proposed approach is an adaptation to 3D of the mesh generation algorithm proposed in \cite{frittelli2021bulk}.

\section{Numerical example on the unit sphere}
\label{sec:example_elliptic_bs_sphere}
\begin{figure}[ht!]
\begin{center}
\hspace*{-11mm}
\includegraphics[scale=0.4]{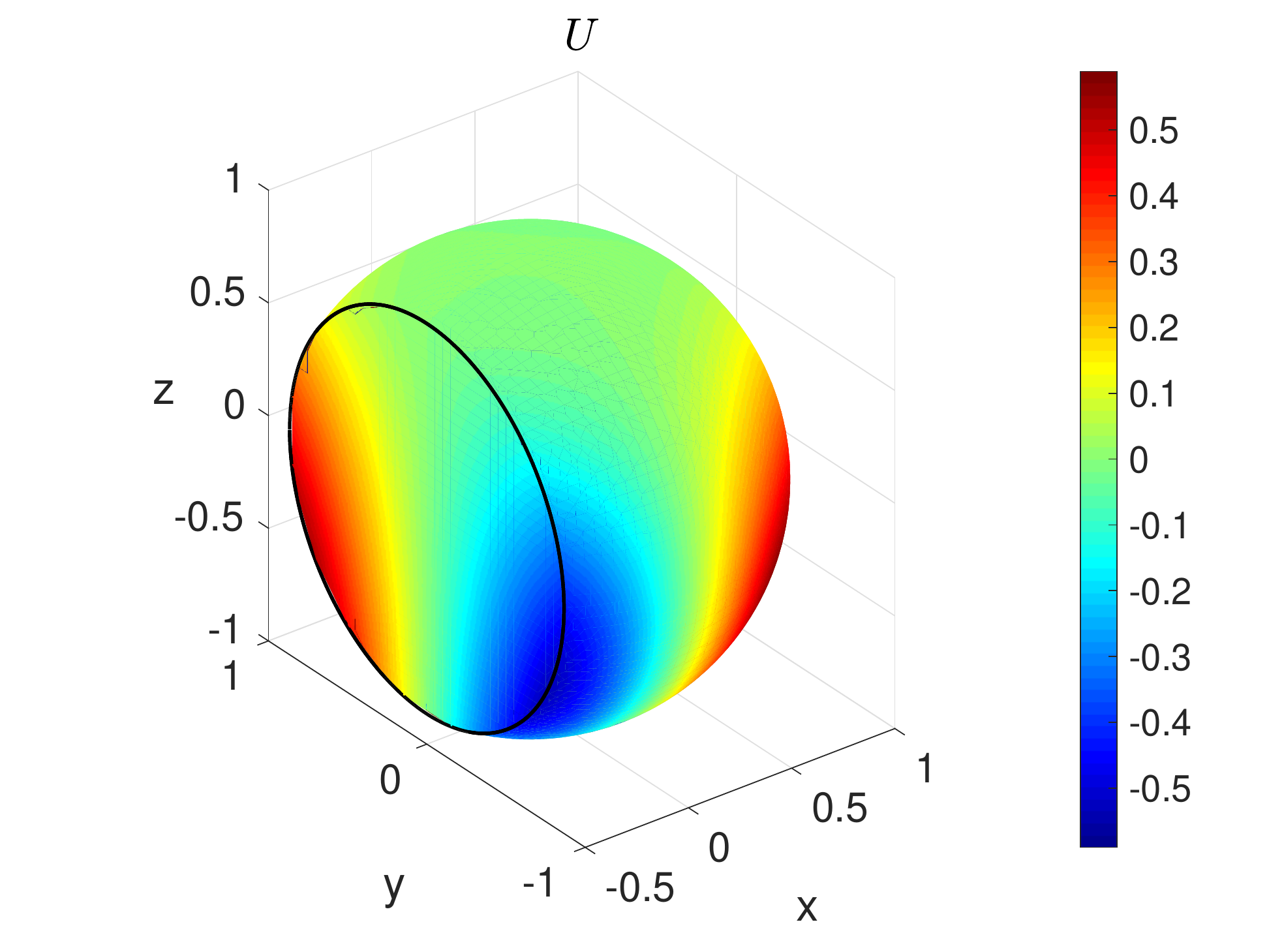}
\hspace*{-5mm}
\includegraphics[scale=0.4]{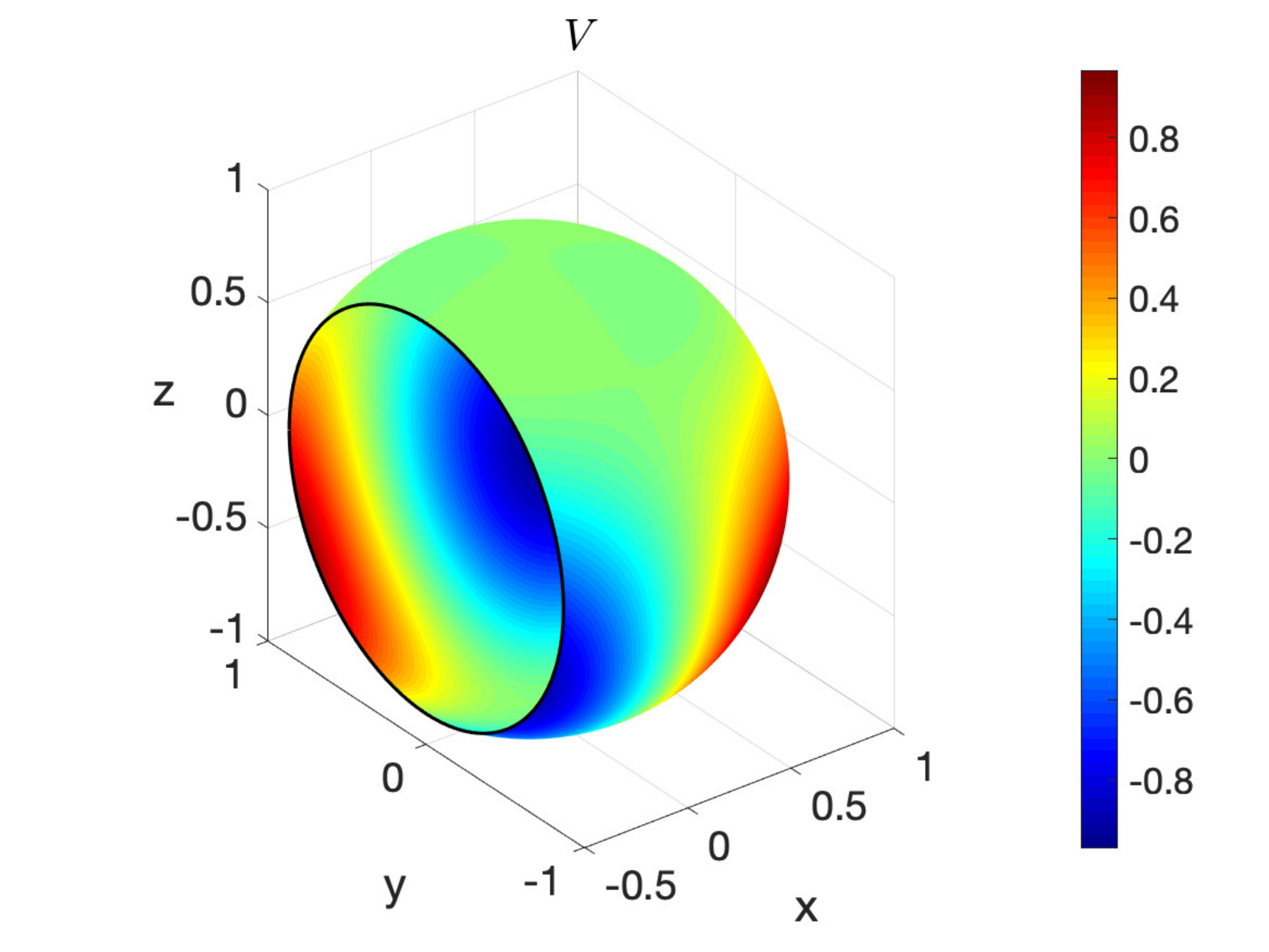}\\
\hspace*{-11mm}
\includegraphics[scale=0.4]{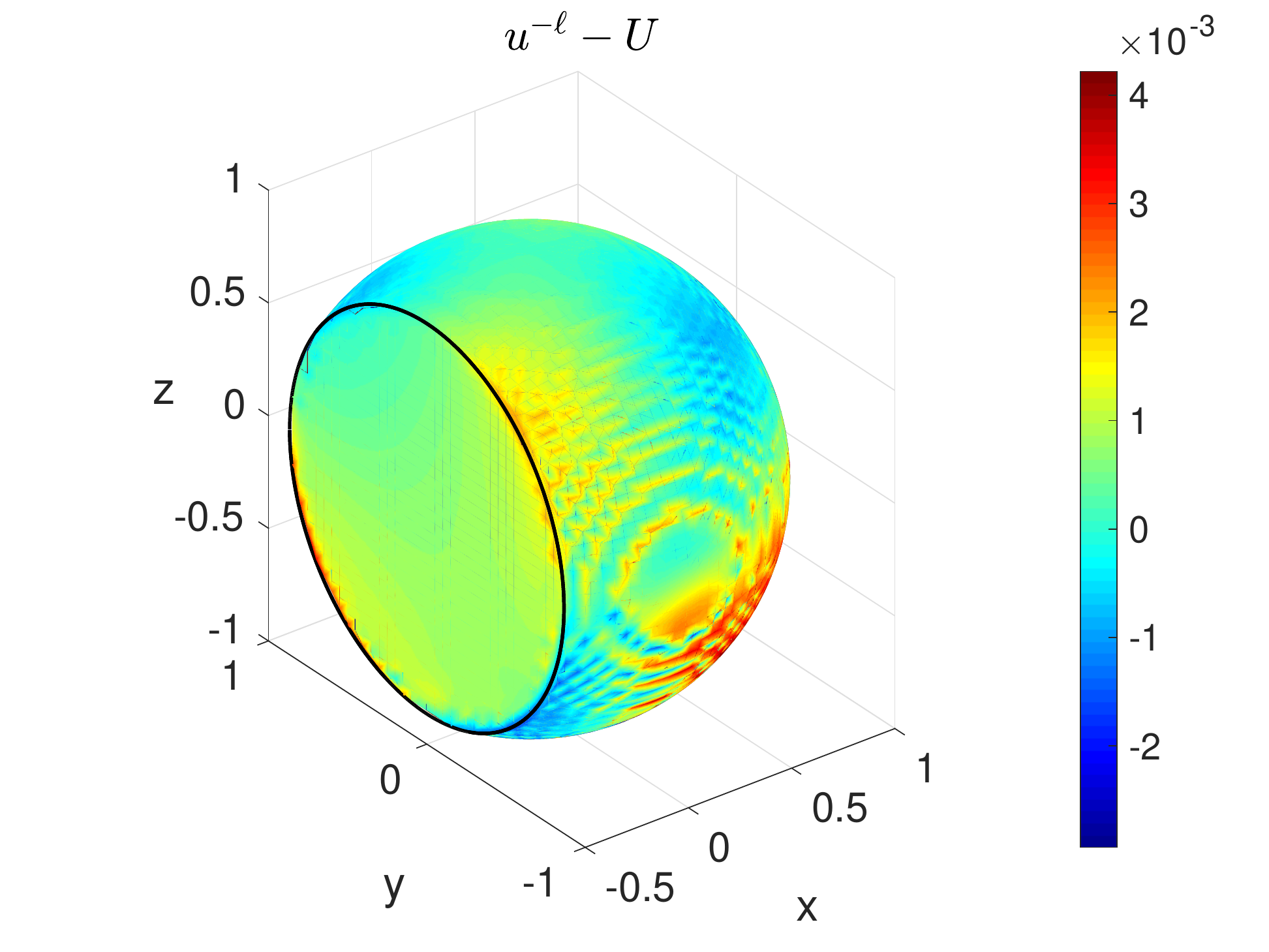}
\hspace*{-5mm}
\includegraphics[scale=0.4]{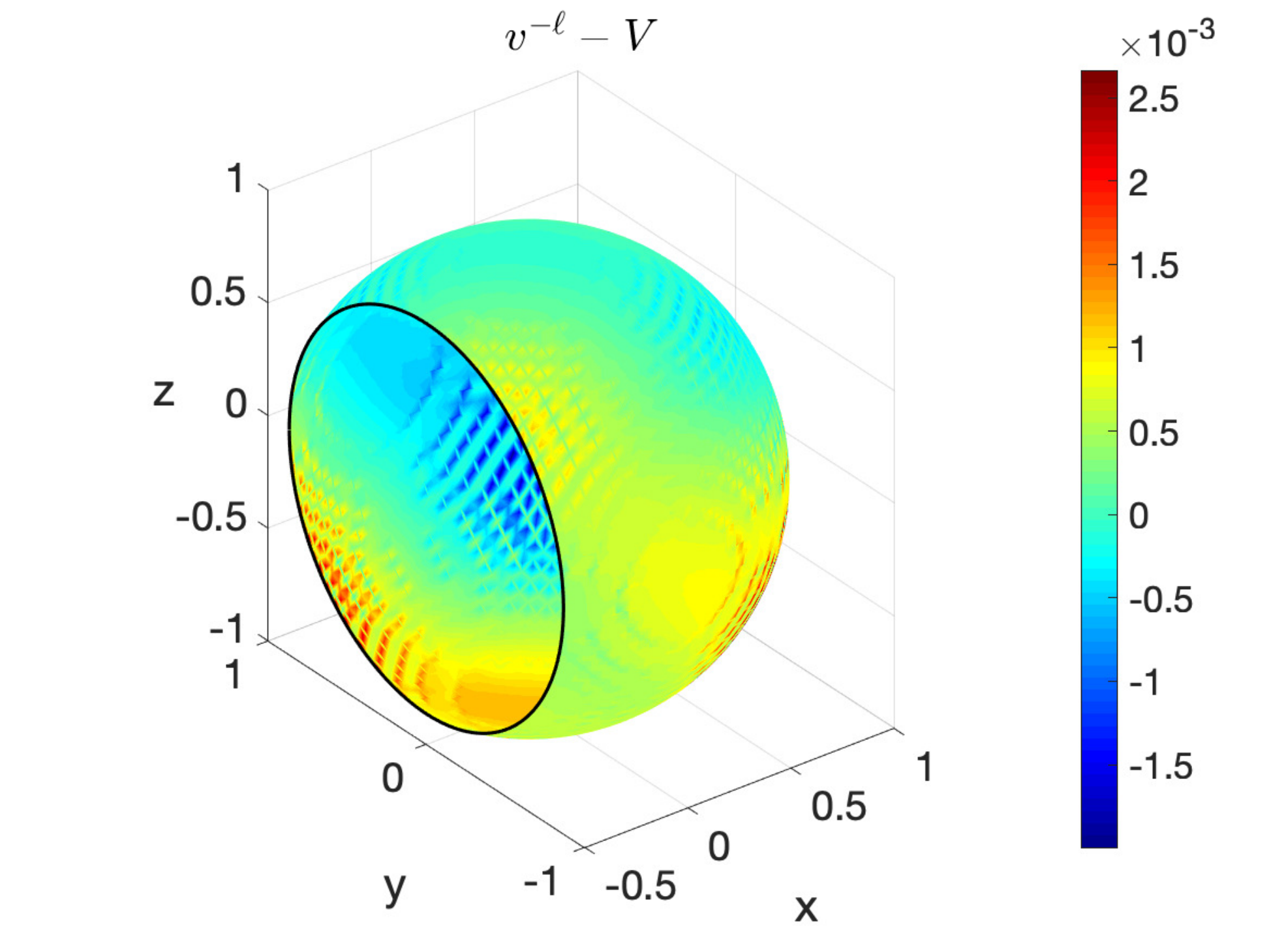}
\end{center}
\caption{Elliptic bulk-surface  problem \eqref{experiment_bs_3d_sphere} on the unit sphere $\Omega$ in 3D: numerical solution obtained on the finest mesh for $i=4$ with $N= 40381$ nodes.  Top row: components $U$ (left) and $V$ (right) of the numerical solution.  Bottom row: pointwise errors in the bulk (left) and on the surface (right).}
\label{fig:bs_3d_numsol_sphere}
\end{figure}

We numerically solve the following elliptic bulk-surface problem on the unit sphere $\Omega$ in 3D:
\begin{equation}
\label{experiment_bs_3d_sphere}
\begin{cases}
-\Delta u + u = xyz - xy\qquad \text{in}\ \Omega;\\
-\Delta_\Gamma v + v +\nabla u\cdot\boldnu = 29xyz - \frac{25}{2}xy\qquad \text{on}\ \partial \Omega;\\
\hspace*{2.5mm} \nabla u\cdot\boldnu = -u + 2v   \qquad \text{on}\ \partial \Omega,
\end{cases}
\end{equation}
whose exact solution is given by $u(x,y,z) = xyz - xy$ for $(x,y,z) \in \Omega$ and $v(x,y,z) = 2xyz - \frac{3}{2}xy$ for $(x,y,z) \in \partial\Omega$.
We consider a sequence of four cubic meshes $i=1,\dots,4$. The $i$-th mesh is obtained by subdividing each dimension into $5i$ intervals, thereby producing a cubic bounding mesh.  From the cubic mesh we obtain a bulk-surface mesh of the sphere as described in Section  \ref{sec:mesh_advantage}.  The coarsest of meshes is shown in Fig.  \ref{fig:mesh_advantage_3}.  On each mesh we solve the discrete problem \eqref{elliptic_problem_BSVEM_linear_system_compact},  we compute the error in $L^2(\Omega)\times L^2(\Gamma)$ norm and the respective convergence rate by the direct solver  \texttt{mldivide} of MATLAB R2019a on a MacBook Pro 2019 with 2,3 GHz 8-Core Intel Core i9 CPU.  As shown in Table  \ref{tab:bs_3d_convergence_sphere}, the convergence in $L^2(\Omega)\times L^2(\Gamma)$ norm is optimal, i.e. quadratic. according to Theorem \ref{thm:BSVEM_L2_convergence}. The numerical solution $(U,V)$ obtained on the finest mesh is plotted in Fig.  \ref{fig:bs_3d_numsol_sphere}, where the bulk component $U$ and the surface component $V$ are shown in separate plots, both cut to show the inside.

\begin{table}
\caption{Elliptic BSPDE \eqref{experiment_bs_3d_sphere} on the unit sphere $\Omega$ in 3D. The BSVEM shows optimal quadratic $L^2$ convergence.  The computational times are shown.}
\begin{center}
\begin{tabular}{@{}ccccccc@{}}
\hline
$i$ & $N$ & $M$ & $h$ & $L^2(\Omega)\times L^2(\Gamma)$ error & EOC & Time (s)\\
\hline
1 & 111 & 56 & 0.6928 &   3.3549e-01 &  -   & 0.002159\\
2 & 799 & 314 & 0.3464 & 5.7422e-02 & 2.5466    & 0.015645\\
3 & 5749 & 1610 & 0.1732 & 1.2235e-02 & 2.2306  & 0.197641  \\
4 & 40381 & 7010 & 0.0866 &  2.8896e-03 & 2.0821 & 5.994934\\
\hline
\end{tabular}
\end{center}
\label{tab:bs_3d_convergence_sphere}
\end{table}

\section{Conclusions}
We have considered a bulk-surface virtual element method (BSVEM) for the numerical approximation of elliptic coupled bulk-surface PDE problems on smooth domains. The proposed method combines a 3D virtual element method (VEM) for the bulk equations (\cite{da2017high}) with a surface virtual element method (SVEM) for the surface equations (\cite{frittelli2018virtual}) and encompasses, in the special case of simplicial bulk-surface meshes,  the BSFEM for bulk-surface RDSs (see e.g. \cite{Madzvamuse_2016}). 

We have introduced polyhedral bulk-surface meshes in three space dimensions and, under minimal mesh regularity assumptions, we have estimated the geometric error arising from domain approximation.  The lack of smoothness in the mapping between the discrete and exact geometries requires the lifting operator to be replaced, in some parts of the analysis, by the Sobolev extension operator.  

The main theoretical result is optimal second-order convergence of the proposed method, provided the exact solution is $H^{2+3/4}$ in the bulk and $H^2$ on the surface.  A relevant by-product is that the lowest order bulk-VEM (\cite{beirao2013basic}) retains optimal convergence even in the simultaneous presence of curved boundaries and non-zero boundary conditions, a result that was not fully addressed in the literature.  The convergence is illustrated with a numerical example on the unit sphere.

We have shown that suitable polyhedral meshes reduce the computational time of mesh generation and matrix assembly from $\mathcal{O}(h^{-3})$ to $\mathcal{O}(h^{-2})$, where $h$ is the meshsize.  This is particularly useful when matrix assembly takes the vast majority of the computational time, i.e.  for (i) time-independent problems and (ii) time-dependent problems on evolving domains, where the matrices must be computed at each timestep.  Polyhedral meshes also allow for simple and efficient adaptive refinement or mesh pasting strategies that would be impossible with tetrahedral meshes,  see for instance \cite{cangiani2014adaptive}. These aspects will be addressed in future studies.

\section*{Acknowledgements}

The work of MF was funded by Regione Puglia (Italy) through the research programme REFIN-Research for Innovation (protocol code 901D2CAA, project number UNISAL026) and by the Italian National Group of Scientific Computing (GNCS-INdAM).
This work (AM) was partly supported by the Global Challenges Research Fund through the Engineering and Physical Sciences Research Council grant number EP/T00410X/1: UK-Africa Postgraduate Advanced Study Institute in Mathematical Sciences, 
the Health Foundation (1902431), the NIHR (NIHR133761) and by an individual grant from the Dr Perry James (Jim) Browne Research Centre on Mathematics and its Applications (University of Sussex).  AM is a Royal Society Wolfson Research Merit Award Holder funded generously by the Wolfson Foundation. AM is a Distinguished Visiting Scholar to the Department of Mathematics, University of Johannesburg, South Africa.
IS is member of the INdAM-GNCS activity group and acknowledges
the PRIN 2017 research Project (No. 2017KL4EF3) ``Mathematics of active materials: from
mechanobiology to smart devices''.

\section*{Conflict of interest}
The authors declare that they have no conflict of interest.

\section*{Data availability}
All data are incorporated into the article.

\bibliographystyle{plainurl}
\bibliography{bibliography}

\appendix
\section{Preliminary definitions and results}

In this Appendix we provide preliminary definitions, results and notations adopted throughout the article. Unless explicitly stated, definitions and results are taken from \cite{dziuk2013finite}.

\subsection{Surfaces and differential operators on surfaces}
Let $\Omega \subset \mathbb{R}^3$ be a compact set such that its boundary $\Gamma := \partial \Omega \subset \mathbb{R}^3$ is a $\mathcal{C}^k$, $k \geq 2$ surface. Since $\Gamma$ can be seen as the zero level set of the \emph{oriented distance function} $d:\mathbb{R}^3 \rightarrow\mathbb{R}$ defined by
\begin{equation*}
d(\boldx) := 
\begin{cases}
-\inf \{\|\boldx-\boldy\| : \boldy\in \Gamma\} \qquad \text{if } \boldx \in \Omega;\\
\hspace*{3mm} 0 \hspace*{34mm} \text{if } \boldx \in\Gamma;\\
\hspace*{3mm} \inf \{\|\boldx-\boldy\| : \boldy \in \Gamma\} \qquad \text{if } \boldx \in \mathbb{R}^3 \setminus \Omega,
\end{cases}
\end{equation*}
then the outward unit vector field $\boldnu:\Gamma \rightarrow\mathbb{R}^3$ can be defined by
\begin{equation}
\label{normal_vector_field}
\boldnu(\boldx) := \frac{\nabla d(\boldx)}{\|\nabla d(\boldx)\|}, \qquad \boldx\in\Gamma.
\end{equation}

\begin{lemma}[Fermi coordinates (\cite{dziuk2013finite})]
\label{lmm:fermi}
If $\Gamma$ is a $\mathcal{C}^k$, $k\geq 2$ surface, there exists an open neighbourhood $U \subset \mathbb{R}^3$ of $\Gamma$ such that every $\boldx\in U$ admits a unique decomposition of the form $\boldx = \bolda(\boldx) + d(\boldx)\boldnu(\bolda(\boldx))$, $\bolda(\boldx) \in \Gamma$. The maximal open set $U$ with this property is called the \emph{Fermi stripe} of $\Gamma$ (see Fig.  \ref{fig:exact_domain}), $\bolda(\boldx)$ is called the \emph{normal projection} onto $\Gamma$ and $(\bolda(\boldx), d(\boldx))$ are called the \emph{Fermi coordinates} of $\boldx$. The oriented distance function fulfils $d\in \mathcal{C}^k(U)$.
\end{lemma}

\begin{definition}[$\mathcal{C}^1(\Gamma)$ functions]
A function $u:\Gamma \rightarrow\mathbb{R}$ is said to be $\mathcal{C}^1(\Gamma)$ if there exist an open neighbourhood $U$ of $\Gamma$ and a $\mathcal{C}^1$ function $\hat{u}:U\rightarrow\mathbb{R}$ such that $\hat{u}_{|\Gamma}  = u$, i.e. $\hat{u}$ is a $\mathcal{C}^1$ extension of $u$ off $\Gamma$.
\end{definition}

\begin{definition}[Tangential gradient and tangential derivatives]
The \emph{tangential gradient} $\nabla_\Gamma u$ of a function $u\in\mathcal{C}^1(\Gamma)$ is defined by $\nabla_\Gamma u(\boldx) := \nabla \hat{u}(\boldx) - (\nabla \hat{u}(\boldx) \cdot \boldnu(\boldx)) \boldnu(\boldx)$ for all $\boldx \in \Gamma$. The result of the computation of $\nablagamma u$ is independent of the choice of the extension $\hat{u}$. The components $D_x u$, $D_y u$ and $D_z u$ of the tangential gradient $\nabla_\Gamma u$ are called the \emph{tangential derivatives} of $u$.
\end{definition}

\begin{definition}[$\mathcal{C}^k(\Gamma)$ functions]
For $k\in\mathbb{N}$, $k>1$, a function $u:\Gamma \rightarrow\mathbb{R}$ is said to be $\mathcal{C}^k(\Gamma)$ if it is $\mathcal{C}^1(\Gamma)$ and its tangential derivatives are $\mathcal{C}^{k-1}(\Gamma)$.
\end{definition}

\begin{definition}[Laplace-Beltrami operator]
The \emph{Laplace-Beltrami} $\Delta_\Gamma u$ of a function $u \in \mathcal{C}^2(\Gamma)$ is defined by $\Delta_\Gamma u(\boldx) := D_xD_x u(\boldx) + D_yD_y u(\boldx)$ for all $\boldx\in\Gamma$.
\end{definition}

\subsection{Bulk- and surface function spaces}
Throughout the paper we will adopt the following notations. For $p \in [1,+\infty]$, $L^p(\Omega)$ and $L^p(\Gamma)$ denote the usual Lebesgue spaces on $\Omega$ and $\Gamma$, respectively, with $\|\cdot\|_{L^p(\Omega)}$ and $\|\cdot\|_{L^p(\Gamma)}$ being the respective norms. For $m \in (0,+\infty)$ and $p\in [1,+\infty]$, $W^{m,p}(\Omega)$ and $W^{m,p}(\Gamma)$ denote the (possibly fractional) Sobolev spaces of order $m$ on $\Omega$ and $\Gamma$, respectively, with $\|\cdot\|_{W^{m,p}(\Omega)}$ and $\|\cdot\|_{W^{m,p}(\Gamma)}$ being the respective norms.  Full definitions can be found in \cite{frittelli2021bulk}.

\begin{lemma}[Inclusion between fractional Sobolev spaces (\cite{Di_Nezza_2012})]
Let $\Omega \subset\mathbb{R}^3$ be a bounded domain with a $\mathcal{C}^1$ boundary $\Gamma$, let $p \in [1,+\infty)$ and $s,s' \in [0,+\infty)$ such that $s < s'$. Then there exists a constant $C>0$ depending on $\Omega$ and $s$ such that
\begin{equation}
\|u\|_{W^{s,p}(\Omega)} \leq C \|u\|_{W^{s',p}(\Omega)},
\end{equation}
for all $u\in W^{s',p}(\Omega)$. Hence, $W^{s,p}(\Omega) \subset W^{s',p}(\Omega)$.
\end{lemma}

\subsection{Fundamental results in bulk- and surface calculus}

\begin{theorem}[Narrow band trace inequality (\cite{elliottranner2013finite})]
With the notations of the previous theorem, there exists $C>0$ depending on $\Omega$ such that any $u \in H^1(\Omega)$ fulfils
\begin{equation}
\label{narrow_band_inequality}
\|u\|_{L^2(U_\varepsilon)} \leq C\varepsilon^{\frac{1}{2}} \|u\|_{H^1(\Omega)}.
\end{equation}
\end{theorem}

\begin{theorem}[Trace theorem and inverse trace theorem (\cite{sobolev1964, Stein_1971})]
\label{thm:trace_theorem}
Let $k\in\mathbb{N}$, $\frac{1}{2} < s \leq k$ and assume that the boundary $\Gamma$ is a $\mathcal{C}^{k}$ surface.\footnote{It is sufficient that $\Gamma$ be a $\mathcal{C}^{k-1,1}$ surface, meaning that its derivatives up to order $k-1$ are Lipschitz continuous. For simplicity, we use the stronger assumption that $\Gamma \in\mathcal{C}^k$.} Then there exists a bounded operator $\Tr: H^s(\Omega) \rightarrow H^{s-\frac{1}{2}}(\Gamma)$, called the \emph{trace operator}, such that $\Tr(u) = u_{|\Gamma}$ and
\begin{equation}
\label{trace_inequality}
\|\Tr(u)\|_{H^{s-\frac{1}{2}}(\Gamma)} \leq C\|u\|_{H^s(\Omega)}, \qquad \forall\ u\in H^s(\Omega).
\end{equation}
The trace operator has a continuous inverse operator $\Tr^{-1}:H^{s-\frac{1}{2}}(\Gamma) \rightarrow H^s(\Omega)$ called \emph{Babi\v{c} inverse} such that
\begin{equation}
\label{inverse_trace_inequality}
\|\Tr^{-1}(v)\|_{H^s(\Omega)} \leq C\|v\|_{H^{s-\frac{1}{2}}(\Gamma)}, \qquad \forall\ v\in H^{s-\frac{1}{2}}(\Gamma).
\end{equation}
\end{theorem}

\begin{theorem}[Sobolev extension theorem (\cite{Stein_1971})]
\label{thm:sobolev_extension_theorem}
Assume that $\Omega \subset\mathbb{R}^3$ has a Lipschitz boundary $\Gamma$, let $r\in\mathbb{N}$ and $p\in [1,+\infty]$. Then, for any function $u\in W^{r,p}(\Omega)$, there exists an extension $\tilde{u} \in W^{r,p}(\mathbb{R}^3)$ such that $\tilde{u}_{|\Omega} = u$ and
\begin{equation}
\label{sobolev_extension}
\|\tilde{u}\|_{W^{r,p}(\mathbb{R}^3)} \leq C\|u\|_{W^{r,p}(\Omega)},
\end{equation}
where $C$ depends on $\Omega$ and $r$, but not on $p$.
\end{theorem}

\begin{theorem}[Sobolev embeddings]
\label{thm:sobolev_embedding}
Let $d\in\mathbb{N}$, $d \geq 2$ be a number of dimensions and assume that $\Omega \subset \mathbb{R}^d$ has a Lipschitz boundary.
\begin{itemize}
\item If $0 < \gamma < 1$, then $H^{d/2+\gamma}(\Omega) \hookrightarrow \mathcal{C}^{0,\gamma}(\Omega)$ is a continuous embedding, hence $\|u\|_{\mathcal{C}^{0,\gamma}(\Omega)} \leq C_\gamma\|u\|_{H^{d/2+\gamma}(\Omega)}$. From the definition of the H\"{o}lder space $\mathcal{C}^{0,\gamma}(\Omega)$ we have that
\begin{equation}
\label{sobolev_inequality_holder}
\|u(\boldx) - u(\boldy)\| \leq C_\gamma\|u\|_{H^{d/2+\gamma}(\Omega)}\|\boldx - \boldy\|^\gamma, \qquad \text{a.e.} \ (\boldx,\boldy) \in \Omega \times \Omega.
\end{equation}
\item If $\varepsilon > 0$, then $H^{d/2+\varepsilon}(\Omega) \hookrightarrow \mathcal{C}(\Omega)$ is a continuous embedding.
\end{itemize}
\begin{proof}
See \cite{adams2003sobolev} for the case of integer-order Sobolev spaces and \cite{Di_Nezza_2012} for the fractional case.
\end{proof}
\end{theorem}

\end{document}